%% file: BPMP.tex
\documentclass[ijoo,nonblindrev]{informs-ijoo}
\OneAndAHalfSpacedXI

\usepackage{comment}

\usepackage{natbib}
 \bibpunct[, ]{(}{)}{,}{a}{}{,}%
 \def\bibfont{\small}%
\usepackage{lineno}
\usepackage{url}
\newtheorem{obs}{Observation}
\newcommand{\anubox}[1]{
\begin{center}
\fbox{
\begin{minipage}{5.8in}
#1
\end{minipage}
}
\end{center}
\smallskip
}

\def\M{\hspace*{1.5em}}
\def\MM{\hspace*{3em}}

\TheoremsNumberedThrough     
\ECRepeatTheorems

\EquationsNumberedThrough    

\MANUSCRIPTNO{}
\begin{document}


\RUNAUTHOR{Dong et al.}

\RUNTITLE{The Backhaul Profit Maximization Problem}

\TITLE{The Backhaul Profit Maximization Problem:  Optimization Models and Solution Procedures}

\ARTICLEAUTHORS{%
\AUTHOR{Yuanyuan Dong}
\AFF{Department of Engineering Management, Information \& Systems, Southern Methodist University, Dallas, TX 75205, USA \EMAIL{ydong@smu.edu}} 
\AUTHOR{Yulan Bai}
\AFF{Department of Engineering Management, Information \& Systems, Southern Methodist University, Dallas, TX 75205, USA \EMAIL{yulanb@smu.edu}}
\AUTHOR{Eli V. Olinick}
\AFF{Department of Engineering Management, Information \& Systems, Southern Methodist University, Dallas, TX 75205, USA \EMAIL{olinick@smu.edu}}
\AUTHOR{Andrew Junfang Yu}
\AFF{Department of Industrial and Systems Engineering, The University of Tennessee, Knoxville, TN 37996-2315, USA \EMAIL{ajyu@utk.edu}}
} 

\ABSTRACT{%
We present a compact mixed integer program (MIP) for the backhaul profit maximization problem in which a freight carrier seeks to generate profit from an empty delivery vehicle's backhaul trip from its last scheduled delivery to its depot by allowing it to deviate from the least expensive (or fastest) route to accept delivery requests between various points on the route as allowed by its capacity and required return time. The MIP is inspired by a novel representation of multicommodity flow that significantly reduces the size of the constraint matrix and the linear programming upper bound on optimal profit compared to a formulation based on the classical node-arc representation. This in turn leads to faster solution times when using a state-of-the-art MIP solver. In an empirical study of both formulations, problem instances with ten potential pickup/dropoff locations and up to 73 delivery requests were solved two times faster on average with our formulation while instances with 20 locations and up to 343 delivery requests were solved 7 to 45 times faster.  The largest instances in the study had 50 locations and 2,353 delivery requests; these instances could not be solved with the node-arc-based formulation, but were solved within an average of less than 40 minutes of real time using our compact formulation.  We also present a heuristic algorithm based on our compact formulation that finds near optimal solutions to the 50-location instances within ten minutes of real time.

}%



\maketitle

\input{section-1-intro.tex}

\input{section-2-node-arc.tex}
\input{section-3-triples.tex}
\input{section-4-heuristic.tex}
\input{section-5-results.tex}

\input{section-6-conclusions.tex}

\begin{APPENDICES}
\input{appendix-A-notation.tex}
\input{appendix-B-mulitple-requests.tex}
\input{appendix-C-expanded-triples-model.tex}
\input{appendix-D-proofs.tex}

\end{APPENDICES}

\bibliographystyle{informs2014}
\bibliography{BPMP}

\end{document}

%% file: section-1-intro.tex
\section{Introduction}
\label{sec:1}
In an increasingly competitive industry, some freight carriers seek to generate profit from an empty delivery vehicle's {\em backhaul} trip from its last scheduled delivery to its depot by allowing it to deviate from the fastest route to accept delivery requests between various points on the route as allowed by its capacity and required return time. This practice has been employed by third party logistics providers (3PLs) in China \citep{dong06}, and is similar to tramp shipping in the maritime cargo transportation industry (e.g., \citealt{CFR}, \citealt{Tramp}) in which ``ships follow the available cargoes, like a taxi" \citep{Bronmo}.
Leading international 3PLs compete on the ability to offer software solutions that allow carriers to solve this type of problem to ``reduce downtime and costly deadhead miles" (\citeauthor{Navisphere} \citeyear{Navisphere,NavisphereVideo}).
Motivated by these applications, we study the {\em backhaul profit maximization problem} (BPMP).

Solving an instance of the BPMP requires simultaneously solving two decision problems: (1) determining a route that the vehicle can take to get from its current location to its depot by the scheduled arrival time, and (2) selecting a profit-maximizing subset of offered point-to-point delivery requests between various points on the route subject to the vehicle's carrying capacity. That is, the vehicle can earn revenue for the 3PL by accepting  delivery requests on the backhaul trip, however the cost incurred for doing so must be taken into consideration. We focus on the case where the vehicle's starting and depot locations are different, as in \cite{dong06}. However, the problem, and our solution approach, can easily extend to the case where they are the same.

There is a rich literature on two classes of problems that are related to the BPMP: pickup and delivery problems (PDP) and vehicle routing problems with profits (VRPwP). In PDP's a fleet of vehicles each starting from, and returning to, the same depot, must fulfill all pickup and delivery requests at minimum total cost. \citeauthor{Berbeglia2007} (\citeyear{Berbeglia2007, Berbeglia20108}) propose an extensive framework for categorizing PDP variants.  However, BPMP does not fall into the framework since all of the delivery requests are optional and the objective in BPMP is to maximize profit rather than minimize the cost of satisfying demand. Additionally, the delivery requests in almost all PDP applications either originate or terminate at the depot.  Another distinguishing feature of the BPMP is that most of the delivery requests originate at points other than the vehicle's starting location and terminate at points other than the depot.

\cite{VRPwProfit} give a comprehensive survey of VRPwP's; a prototypical example application is a variation of the travelling salesman problem (TSP) in which the salesman receives varying amounts of revenue for visiting each city in the problem instance and is not required to visit all of the cities. Thus, the objective is to find a tour that maximizes the profit resulting from the revenue obtained by visiting a subset of the cities minus the travel cost. The motivation for a VRPwP  can be similar to  that of BPMP; for example, the revenue from visiting a city might come from making a delivery from a depot. However, most of the VRPwP variants discussed in \cite{VRPwProfit} are structurally distinct from BPMP in part because the BPMP considers a variable travel cost based on the vehicle's load in addition to the fixed city-to-city travel cost in TSP-type applications. Although the PDP, VRPwP, and BPMP are different problem classes, they all contain the well-studied orienteering problem (OP) as a special case.

The OP is motivated by a type of cross-country running race in which competitors earn points for visiting designated locations between the start and finish lines, and has received considerable attention in the literature.  Typically there is a time limit for the race, and so the problem is to find a route from the start to the finish that maximizes total points earned within the allowable time \citep{Tsiligirides}. For recent, comprehensive surveys on the OP and its many variants see  \cite{Vansteenwegen20111} and  \cite{Gunawan2016315}. In Section \ref{sec:complexity} we use the fact that the OP is ${\cal NP}$-hard to establish the computational complexity of the BPMP.

To the best of our knowledge, \cite{Yu} presented the first exact solution procedures for BPMP: a mixed integer programming (MIP) model and an exhaustive search procedure. The MIP is based on a node-arc formulation of multicommodity flow and requires considerable computing resources (time and memory) to solve. For example, the MIP for a 30-location instance in \citep{Yu} has over 600,000 variables and 27,000 constraints.  In this paper we present a new MIP formulation of the BPMP based on a compact representation of multicommodity flow proposed by  \cite{NET:NET21583}. Our new model for the BPMP significantly reduces the number of constraints and binary variables. Furthermore, we demonstrate empirically that the new formulation has a much stronger linear programming relaxation, and present computational results demonstrating that CPLEX can solve the new model significantly faster than the node-arc model proposed by \cite{Yu}. We also present a straight-forward heuristic based on our new formulation that is found experimentally to be 1.79 to 5.75 times faster for 40-node instances and 2.16 to 14.79 times faster for  50-node instances than the exact method with minimial loss in solution quality.

The rest of this paper is structured as follows. In Section \ref{sec:BPMP}, we review the MIP model of the BPMP proposed by \cite{Yu}, show how the MIP can be improved (e.g. with valid inequalities and preprocessing), and prove that the BPMP is ${\cal NP}$-hard.
In Section, \ref{sec:triples} we present our new MIP formulation, which we call the {\em triples formulation}, and show that it is significantly more compact than the node-arc formulation. We propose our heuristic for the BPMP based on the triples formulation in Section \ref{sec:heuristic}. In Section \ref{sec:results}, we summarize an extensive computational study solving instances of the BPMP with the node-arc and triples formulations, and the heuristic. As we note in our conclusions in Section \ref{sec:conclusions}, the results in Section \ref{sec:results} suggest that the compact representation of multicommodity flow exploited by the triples formulation has the potential to significantly improve the efficiency of solving other logistics problems related to the BPMP.

%% file: section-2-node-arc.tex
\section{The BPMP}
\label{sec:BPMP}

Figure 1 gives a graphical illustration of a BPMP instance.  In the figure, an empty vehicle is scheduled to travel from location 1 to location 10, and the eight other locations in the figure are potential stops to pickup and/or deliver cargo. The numbers next to the arcs represent the weight of the delivery requests in tons, and the vehicle's capacity is one ton. The feasible solution shown in the figure routes the vehicle on the indicated route, $1 \rightarrow 6 \rightarrow 8 \rightarrow 10$, and accepts four delivery requests: (1) 0.8 tons from location 1 to location 6, (2) 0.2 tons from location 1 to location 10, (3) 0.6 tons from location 6 to location 10, and (4) 0.2 tons from 8 to 10.  In this section we give a formal definition of the BPMP and review the node-arc formulation proposed in \cite{Yu}. We then show how the node-arc formulation can be strengthened,  and establish the computational complexity of the BPMP.

\begin{figure}[!h]
    \begin{center}
   \includegraphics[scale=.6]{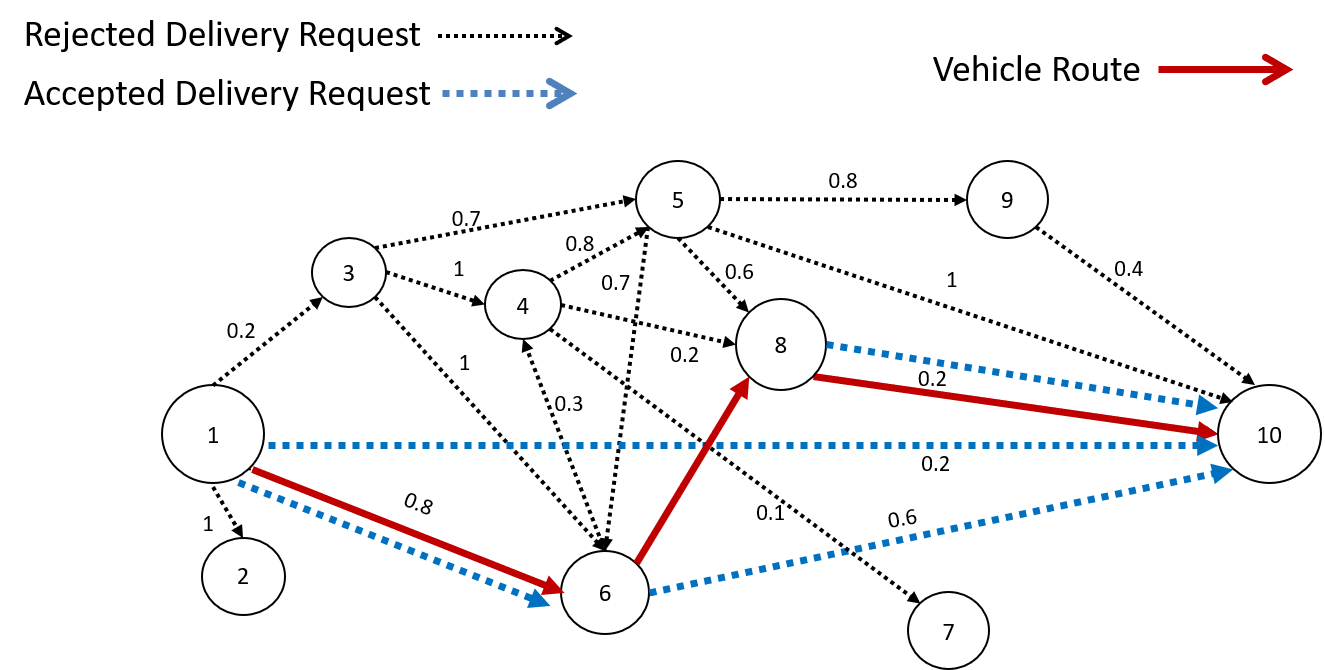}
   \end{center}
  \caption{BPMP Example}
  \label{fig:1}
\end{figure}

\subsection{Formal Problem Definition}
\label{sec:2.1}

An instance of the BPMP is defined on a network with node set $\mathcal{N} = \{1, 2, \ldots, n\}$ and arc set $\mathcal{A}= \{(i, j) | i < n, j > 1, i \neq j \}$. Node 1 represents the vehicle's current location, node $n$ represents the depot, and nodes $2, 3, \ldots, n-1$ represent potential customer locations for optional pickup and delivery on the way to the depot. Hereinafter, we use the terms node and location interchangeably. Each arc $(i, j)$ has a nonnegative distance $d_{ij}$ representing the driving distance, in miles, from node $i$ to node $j$. The arc distances are assumed to be Euclidean (or at least to satisfy the triangle inequality) as is common in the logistics literature \citep{Berbeglia2007}. The BPMP proposed by \cite{Yu} is a generalization of the problem described by \cite{dong06} in which all of the delivery requests are for full truckloads.  \cite{Yu} consider requests that are less than the vehicle's capacity, but assume that there is at most one delivery request per pair of pickup and dropoff locations as in \cite{dong06}.  Hence, we denote the set of delivery requests as $\mathcal{R} \subseteq A$ and the weight, in tons, of a load that a potential customer would like to ship from node $i$ to node $j$ as $w_{ij}$. In Appendix \ref{appendix-multiple-requests}, we describe a straight-forward adaptation of the network so that the model may be applied to instances where there are multiple requests between particular node pairs.  If the vehicle accepts the delivery request $(i, j) \in \mathcal{R}$, then the customer pays $p$ $d_{ij}$  $w_{ij}$, where $p$ is the price charged (revenue received) in dollars per mile per ton. Thus, customers are not charged for any detours the vehicle makes between the end nodes of their delivery requests. In Figure \ref{fig:1}, for example, the customer shipping from location 1 to location 10 is charged $p$ $d_{1,10}$ $w_{1, 10}$ even though their cargo travels a distance of $d_{16} + d_{68} + d_{8,10}$. For convenience in writing the formulations, we let $w_{ij} = 0$ for $(i,j) \in \mathcal{A} \setminus \mathcal{R}$.

The vehicle weighs $v$ tons when empty and has a carrying capacity of $Q$ tons. The vehicle incurs a travel cost (fuel, wear and tear, etc.) of $c$ dollars per mile per ton. Thus, an empty vehicle incurs a cost of $c$ $v$ $d_{ij}$ when traversing arc $(i, j)$,  while a fully loaded vehicle incurs a cost of $c$ $(v + Q)$ $d_{ij}$. The vehicle has $\tau$ hours to reach node $n$ given the deadline for arriving at the depot; assuming that the vehicle travels at a known average speed, the deadline is enforced by limiting the  distance the vehicle travels on its route from node 1 to node $n$ to at most $D$ miles.

It is assumed that all delivery requests are available on a ``spot" market for pickup and delivery within the next $\tau$ hours and that partial delivery is not allowed for any request. We also assume that the cargoes are  heterogeneous and may all be transported together in the vehicle. Thus, weight is the only factor considered in the capacity constraint. The objective of the problem is to determine a selection of delivery requests to accept and a corresponding route that maximizes the total profit subject to the distance  and capacity limits.

\subsection{Node-Arc Formulation of BPMP}
\label{sec:node-arc}

In this section, we review the node-arc formulation of BPMP proposed by \cite{Yu}. The formulation uses the binary variable $x_{ij}$ to indicate whether or not the vehicle traverses arc $(i, j)$, and binary variable $y_{kl}$ to indicate whether or not to accept delivery request $(k, l)$. Binary variable $z_{kl,ij}$ determines whether or not delivery request $(k, l)$ is performed via arc $(i, j)$. In multicommodity flow terms, each delivery request is a commodity and variable $z_{kl,ij}$ indicates if the commodity shipped from node $k$ to node $l$ flows on arc $(i, j)$.  Variable $\theta_{ij}$ represents the total flow (i.e., tons of cargo) transported on arc $(i, j)$. Sequence variables $s_i \ge 0$, for $i= 1, \ldots, n$, track the relative order in which nodes are visited. Using the notation above, which is summarized in Appendix \ref{appendix-notation}, the node-arc formulation for BPMP maximizes profit given by

\begin{equation}
 p \sum_{(k, \ell) \in \mathcal{R}} d_{k\ell} w_{k\ell} y_{k\ell} - c \sum_{(i, j) \in \mathcal{A}}  d_{ij} \theta_{ij} - c v \sum_{(i, j) \in \mathcal{A}} d_{ij} x_{ij}
\label{eqn:1}
\end{equation}

\noindent subject to

\begin{eqnarray}
  \sum_{j=2}^{n} x_{1j} & = & 1 \label{eqn:2e} \\
  \sum_{i = 1}^{n-1} x_{in} & = & 1 \label{eqn:2f} \\
  \sum_{i \in \mathcal{N} \setminus \{k,n\}} x_{ik} & = & \sum_{j \in \mathcal{N} \setminus \{1, k\}} x_{kj} \qquad \forall k \in \mathcal{N} \setminus \{1, n\} \label{eqn:2g}\\
  \sum_{(i,j) \in \mathcal{A}} d_{ij} x_{ij} & \le & D \label{eqn:2h} \\
  \sum_{i \in \mathcal{N} \setminus \{k, n\}} x_{ik} & \le & 1 \qquad \forall k \in \mathcal{N} \setminus \{1, n\} \label{eqn:node-degree} \\
  s_i - s_j + (n + 1) x_{ij} & \le & n   \qquad \forall (i, j) \in \mathcal{A} \label{eqn:2k} \\ \nonumber \\
  \sum_{(k,l) \in \mathcal{R}} z_{kl,ij} & \le & M x_{ij} \qquad \forall (i, j) \in \mathcal{A} \label{eqn:2d} \\
  \sum_{j \in \mathcal{N} \setminus \{1, k\}} z_{kl,kj} &  = &  y_{kl}   \qquad \forall (k, l) \in \mathcal{R} \label{eqn:2a} \\
  \sum_{i \in \mathcal{N} \setminus \{l,n\}} z_{kl,il} &  = &  y_{kl} \qquad \forall (k, l) \in \mathcal{R} \label{eqn:2b} \\
 \sum_{\{i \in \mathcal{N}: (i, h) \in \mathcal{A}\}} z_{kl,ih} &  = & \sum_{\{j \in \mathcal{N}: (h, j) \in \mathcal{A}\}} z_{kl,hj}  \qquad  \forall (k, l) \in \mathcal{R}, h \in \mathcal{N} \setminus \{k,l\} \label{eqn:2c} \\
  \theta_{ij} & = & \sum_{(k,l) \in \mathcal{R}} w_{kl} z_{kl,ij} \qquad (i, j) \in \mathcal{A} \label{eqn:2i} \\ \nonumber \\
  \theta_{ij} & \le & Q \qquad \qquad \forall (i, j) \in \mathcal{A} \label{eqn:2j} \\
      s_i   & \ge 0 & \qquad \qquad  \forall i \in {N} \label{eqn:s_ge_0}\\
  x_{ij} & \in & \{0, 1\} \qquad \forall (i, j) \in \mathcal{A} \label{eqn:binary_x} \\
  y_{ij} & \in & \{0, 1\} \qquad \forall (i, j) \in \mathcal{R} \label{eqn:binary_y} \\
  z_{kl,ij} & \in & \{0, 1\} \qquad \forall (k,l) \in \mathcal{R}, (i, j) \in \mathcal{A} \label{eqn:binary_z}
\end{eqnarray}

The three terms in the objective function (\ref{eqn:1}) are revenue from accepted delivery requests,

\[p \sum_{(k, l) \in \mathcal{R}} d_{kl} w_{k l} y_{k l},\]

\noindent travel cost related to delivery requests (i.e., cargo-carrying costs),

\[c \sum_{(i, j) \in \mathcal{A}} d_{ij} \theta_{ij}, \]

\noindent and the vehicle-related travel cost,

\[c v \sum_{(i, j) \in \mathcal{A}} d_{ij} x_{ij}.\]

Constraints (\ref{eqn:2e})--(\ref{eqn:2k}) are the {\em vehicle-routing constraints}. The vehicle's route is represented as a unit of flow from node 1 to node $n$ by constraints (\ref{eqn:2e})--(\ref{eqn:2g}), and constrained to at most $D$ miles by constraint (\ref{eqn:2h}). The {\em node-degree}  (\ref{eqn:node-degree}), and {\em subtour elimination} constraints (\ref{eqn:2k}) adapted from the sequential formulation for the TSP proposed by \citep{miller} ensure that the vehicle follows a simple path from node 1 to node $n$. The sequence variables determine the relative order in which nodes are visited by the vehicle. This approach is used to eliminate subtours instead of the conventional approach used by \cite{Dantzig} in order to reduce problem size \citep{Yu}. Note that in the absence of the node-degree constraints (\ref{eqn:node-degree}), it is possible that the positive $x$ variables in the LP relaxation of the node-arc formulation form a vehicle route with subtours. Thus, constraint set (\ref{eqn:node-degree}) is intended to strengthen the formulation.

The logical connection between $x_{ij}$ and $z_{kl,ij}$ is enforced by constraint set (\ref{eqn:2d}). If the left-hand side is positive, then $x_{ij} = 1$; intuitively, this means that if at least one of the accepted delivery requests is routed via arc $(i, j)$ then the vehicle must travel on that arc.  Conversely, if arc $(i, j)$ is not on the vehicle's route (i.e., $x_{ij} = 0$), then $\sum_{(k, l) \in \mathcal{R}} z_{kl,ij} = 0$.

The next group of constraints, (\ref{eqn:2a})--(\ref{eqn:2j}), model the movement of cargo carried by the vehicle as multicommodity flow. Constraint sets (\ref{eqn:2a}) and (\ref{eqn:2b}) enforce the logical relationship between $y_{kl}$ and $z_{kl,ij}$;  if delivery request $(k, l)$ is accepted, then the vehicle's route must contain an arc leaving node $k$ and an arc entering node $l$.  Constraints  (\ref{eqn:2c}) are flow-conservation constraints for intermediate nodes on the path the vehicle takes from node $k$ to node $l$. Together with constraints (\ref{eqn:2e})--(\ref{eqn:2g}), these constraints ensure that if delivery request $(k, l)$ is accepted, then the vehicle must visit node $k$ before visiting node $l$. Constraint set (\ref{eqn:2i}) determines the total load of the vehicle (tons carried) on each arc and ensures that the arc flows are nonnegative. The capacity limit is enforced by constraint set (\ref{eqn:2j}).  We denote a solution to the node-arc formulation by a tuple of unsubscripted variables $(x, y, z, s)$.

\subsection{Strengthening the Node-Arc Formulation}
\label{sec:best_node_arc}

The subtour elimination constraints, (\ref{eqn:2k}), and the constraints linking the $x$ and $z$ variables, (\ref{eqn:2d}), are salient, potential sources of weakness in the node-arc model.  In a preliminary study \citep{BO}, we found that lifting (\ref{eqn:2k}) \'{a} la \cite{DLMTZ} strengthened the node-arc formulation but actually lead to longer solution times in many instances.
\cite{Dong} experimented with replacing constraint (\ref{eqn:2d}) with $z_{kl,ij}  \le x_{ij}$ for all combinations of $(k,l) \in \mathcal{R}$ and $(i, j) \in \mathcal{A}$.  This ``dissagregated" formulation is stronger, but a factor of $n$ larger than the original
node-arc formulation in terms of the number of constraints. \cite{Dong} found that disaggregating constraint (\ref{eqn:2d}) reduced the upper bound on profit from the LP relaxation significantly (up to 93\%), but had a counterproductive effect of increasing solution time in all cases.
\cite{Yu} do not specify the ``big-$M$" value they used for the right-hand side of the constraint (\ref{eqn:2d}).  \cite{BO} derive a data-independent value of $\frac{n^2-n}{2}$, which is the value we used in our preliminary study. For any particular problem instance a tighter bound might be found by solving a binary knapsack problem that maximizes the number of delivery requests accepted subject to the total weight of the accepted requests being at most the vehicle capacity, $Q$.

Testing 20-node instances from \cite{Yu}, we found that CPLEX took between 5 minutes to 35 minutes to solve the node-arc model proposed in \citep{Yu}. In a comprehensive study applying our own insights and adapting techniques from the literature on related problems (e.g., \cite{FGT}), we reduced the solution-time range to 31 to 105 seconds \citep{BO}. Adopting the best practices from our study \citep{BO}, we introduce the {\em enhanced node-arc formulation} as follows.

The original node-arc formulation proposed by \cite{Yu} uses constraint  set (\ref{eqn:2j}), $\theta_{ij} \le Q$, to ensure that the total amount of flow on arc $(i, j)$ does not exceed the vehicle's capacity. Notice, however, that if the vehicle does not travel on arc $(i, j)$, there should be no flow on the arc. If the vehicle does travel on arc $(i, j)$, then the flow on the arc can be at most $Q$. Therefore, constraint (\ref{eqn:2j}) can be replaced by the following {\em conditional arc-flow} constraints:

\begin{equation}
  \theta_{ij}  \le  Q x_{ij} \qquad (i,j) \in {\cal A}. \label{eqn:conditional arc flow}
\end{equation}

\noindent Since constraint (\ref{eqn:2i}) defines the arc flows in terms of the $z$ variables,  adopting the conditional arc-flow constraints makes the constraints  linking the $x$ and $z$ variables, (\ref{eqn:2d}),  redundant. Using the bound from \citep{BO}, $M = 45$, 190, 435, and 780, for $n = 10$, 20, 30, and 40, respectively. In practice, $Q$ is on the order of  $40$ for vehicles traveling on the National Highway System in the United States \citep{WeightLimit}. Thus, replacing constraints (\ref{eqn:2d}) and (\ref{eqn:2j}) with constraint (\ref{eqn:conditional arc flow}) strengthens the formulation, and in  our preliminary study \citep{BO} we found that the resulting model was faster to solve.
In our preliminary experiments \citep{BO} we found that CPLEX's performance improved when we dropped the node-degree constraints (\ref{eqn:node-degree}), which are not are not strictly necessary for the formulation to be valid due to the given the subtour elimination constraints (\ref{eqn:2k}).

Hereinafter, we refer to the MIP resulting from making the following changes to the
node-arc formulation as the {\em enhanced node-arc formulation}: (i)  drop the node-degree constraints (\ref{eqn:node-degree}), and (ii)
replace the $x-z$ linking constraints (\ref{eqn:2d}) and capacity constraints (\ref{eqn:2j})  with the conditional arc-flow constraints (\ref{eqn:conditional arc flow}).

\subsection{Computational Complexity of BPMP}
\label{sec:complexity}

\cite{Yu} discuss the exponential growth of the BPMP solution space as a function of the number of locations in the problem and give an example four-node BPMP instance with locations that has 146 distinct combinations of routes and accepted delivery requests.  Note that this is a larger solution space than a four-city TSP, which has 24 distinct tours. A 10-node TSP has approximately $3.6 \times 10^6$ solutions while a 10-node BPMP has on the order of $10^{18}$ potential solutions. In general, the number of potential solutions for a BPMP instance with $k$ locations between the starting and destination locations is $\sum_{r=0}^{k} P^r_k \times 2^{C^2_{r+2}}$, where $r$ is the number of locations selected for the route, $P^r_k$ is the number of permutations corresponding to routes visiting $r$ out of $k$ locations, and $2^{C^2_{r+2}}$ is the number of different ways of choosing an ordered pair from a set of $r+2$ locations \citep{Yu}. Thus, even relatively small problem instances can be quite challenging.  We conclude this section with a polynomial-time reduction of the orienteering problem described in Section \ref{sec:1} to the BPMP and hence show that the BPMP is an ${\cal NP}$-hard optimization problem.

\begin{theorem}
The BPMP is ${\cal NP}$-hard.
\end{theorem}

\proof{Proof}
It is straight-forward to show that the size of the node-arc formulation of BPMP is bounded by a polynomial function of $n$. Hence, BPMP is in the problem class ${\cal NP}$. An instance of the orienteering problem, which is known to be ${\cal NP}$-hard \citep{Golden},  can be modeled as a special case of the BPMP in which every delivery request is destined for the depot as follows:

\begin{itemize}
\item Let nodes $1$ and $n$ in the BPMP instance represent the starting and finishing lines of the OP instance, respectively, and let nodes $2, 3, \ldots, n-1$ represent the designated points that the racer can visit to earn points.

\item Let $\tau$ be the maximum time allowed for the race (in hours). Without loss of generality assume that the racer runs at an average rate of one mile per hour so that $d_{ij}$ in the BPMP instance is equal to the given number of hours required for the racer to run from location $i$ to location $j$ in the OP instance. Hence, the corresponding BPMP distance limit is $\tau$ miles.

\item Let $P_i$ denote the points earned for visiting location $i$ in the OP instance. In the BPMP instance
let $w_{in} = \frac{P_i}{d_{in}}$  for $1 \le i < n$,  and $w_{ij} = 0$ for all other $(i,j) \in A$.

\item Let the travel cost  $c$ in the BPMP instance equal zero dollars per mile.

\item Let the empty vehicle weight $v$ equal in the BPMP instance equal one ton.

\item Let the price charged per mile per ton in the BPMP instance be $p = 1$ dollar.

\item Since capacity is irrelevant in OP, let $Q = \sum_{1 \le i < n} w_{in}$ in the BPMP instance.

\end{itemize}

\noindent Clearly, a feasible solution to the BPMP instance corresponds to a feasible route for the racer in the OP instance. Observe that the BPMP objective function (\ref{eqn:1}) simplifies to $\sum_{i=1}^{n-1} P_i y_{in}$. Thus, the profit earned by the vehicle in the BPMP instance is equal to the total points awarded for the corresponding route in the OP instance. It follows that the BPMP is ${\cal NP}$-hard.
\Halmos
\endproof 

%% file: section-3-triples.tex
\section{Triples Formulation of the BPMP}
\label{sec:triples}

The MIP in Section \ref{sec:BPMP} is based on the node-arc formulation for multicommodity flow. Our triples formulation of the BPMP is adapted from a compact representation of multicommodity flow that has been successfully applied to the maximum concurrent flow problem \citep{NET:NET21583}. The triples formulation is based on a description of flow where {\em triples variable} $u_{ij}^k$ for node triple $(i,j,k)$ represents the total flow on all paths from node $i$ to node $j$ with arc $(i,k)$ as the first arc.  In the BPMP application, node $k$ in triple $(i,j,k)$ cannot be the starting location (node 1) and cannot be the depot (node $n$). When $u_{ij}^k$ is positive, we say that $u_{ij}^k$ units of flow from node $i$ to node $j$ are {\em diverted} through node $k$. It is important to note that variable $u_{ij}^k$ does not specify how the flow travels from node $k$ to node $j$, and that flow that is not diverted (i.e., {\em direct flow}) is represented implicitly; that is, there is no variable in the triples formulation that explicitly represents the amount of flow from node $i$ to node $j$ sent on arc $(i, j)$.  We begin this section with an illustrative triples solution to an example BPMP instance in Subsection \ref{sec:TriplesExample}. We then present the triples formulation for BPMP and establish its validity in Subsection \ref{sec:TriplesMIP}. In Subsection \ref{sec:enhanced triples} we present the {\em enhanced triples formulation}, which is stronger and even more compact than the initial formulation presented in Subsection \ref{sec:TriplesMIP}.  We conclude this section by using the enhanced triples formulation to derive an upper bound on profit  on the BPMP in Subsection \ref{sec:upper bound on profit}.

\subsection{Representing Flow with Triples Variables}
\label{sec:TriplesExample}

\begin{figure}[ht]
    \begin{center}
   \includegraphics[width=0.6\textwidth]{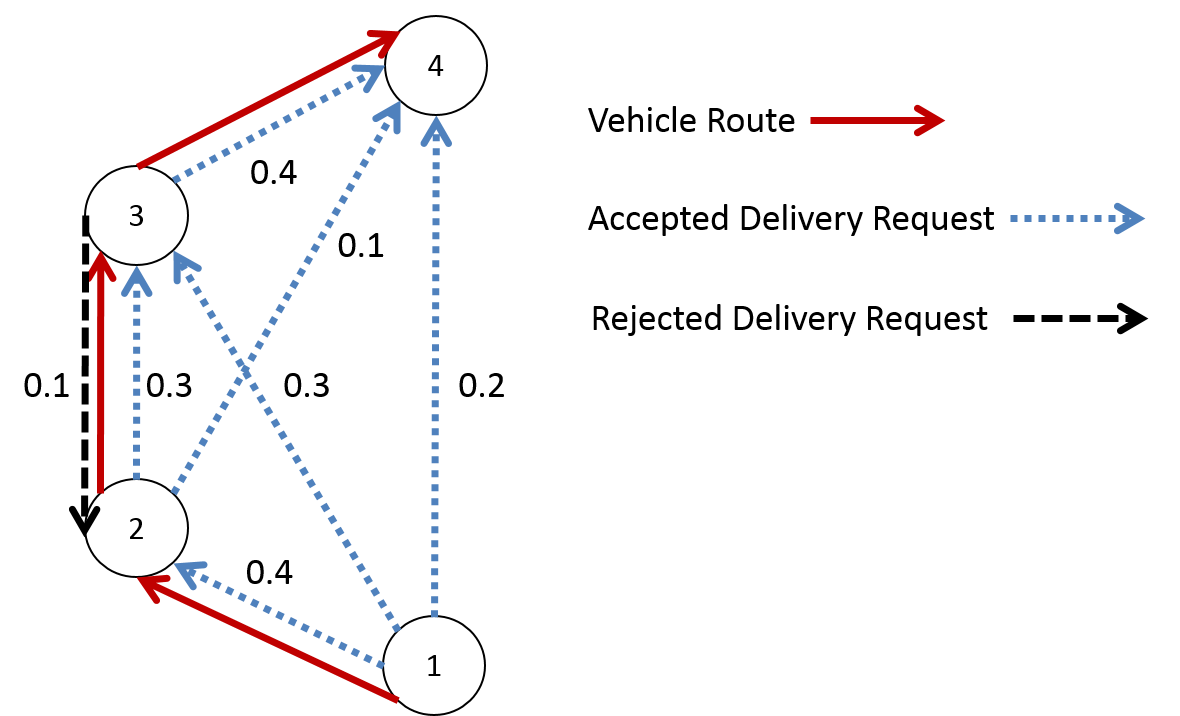}
   \end{center}
  \caption{Example Solution to a BPMP Instance with $n = 4$}
  \label{fig:2}
\end{figure}

To illustrate the representation of the flow of cargo with triples variables, consider the four-node example BPMP instance and solution shown in Figure \ref{fig:2}. As in Figure \ref{fig:1}, the numbers next to the arcs are the weights of the delivery requests (in tons) and the vehicle has a one-ton capacity. As in the node-arc formulation, the vehicle route indicated in Figure \ref{fig:2} is represented with positive $x$ variables, $x_{12} = x_{23} = x_{34} = 1$, and the accepted delivery requests are represented with positive $y$ variables, $y_{12} = y_{13} = y_{14} = y_{23} = y_{24} = y_{34} = 1$.  Accepted delivery requests between consecutive nodes on the vehicle's path, such as the request from node 1 to node 2, do not need to be diverted through intermediate nodes.  Therefore, all triples variables of the form $u_{12}^k$,  $u_{23}^k$, and $u_{34}^k$ are equal to zero.
The diverted flow for the delivery requests $(1, 3)$, $(1, 4)$, and $(2, 4)$ is represented with positive triples variables $u_{13}^2 = 0.3$, $u_{14}^2 = 0.2$, and $u_{24}^3 = 0.3$, which are interpreted as follows:

\begin{itemize}
\item  $u_{13}^2 = 0.3$ indicating that the flow for delivery request $(1, 3)$ is diverted through node 2. Thus, all flow from node 1 to node 3 is sent on arc $(1, 2)$ and is then sent from node 2  to node 3. As noted above, all flow in the solution shown in Figure \ref{fig:2} from node 2 to node 3 is direct flow. Thus, the delivery request from node 1 to node 3 follows the path composed of arcs $(1, 2)$ and $(2, 3)$.

\item $u_{14}^2 = 0.2$ indicating that the flow for delivery request $(1, 4)$ is diverted through node 2. Thus, all flow from node 1 to node 4 is sent on arc $(1, 2)$  and is then sent from node 2 to node 4 via an unspecified path.

\item  $u_{24}^3 = 0.3$ indicating that the solution sends 0.3 units of flow from node 2 to node 4 by diverting it through node 3. Since all flow from node 3 to node 4 is direct, the solution sends 0.3 units of flow on the path composed of arcs $(2, 3)$ and $(3, 4)$. Note that $u_{24}^3$ is the combined flow from two different delivery requests: $(1, 4)$ with $w_{14} = 0.2$ and $(2, 4)$ with $w_{24} = 0.1$.  Thus, the triples variables indicate that the flow for delivery request (1, 4) follows the path $1 \rightarrow 2 \rightarrow 3 \rightarrow 4$, and the flow for delivery request (2, 4) follows the path $2 \rightarrow 3 \rightarrow 4$.

\end{itemize}

\begin{figure}[ht]
    \begin{center}
   \includegraphics[scale=.85]{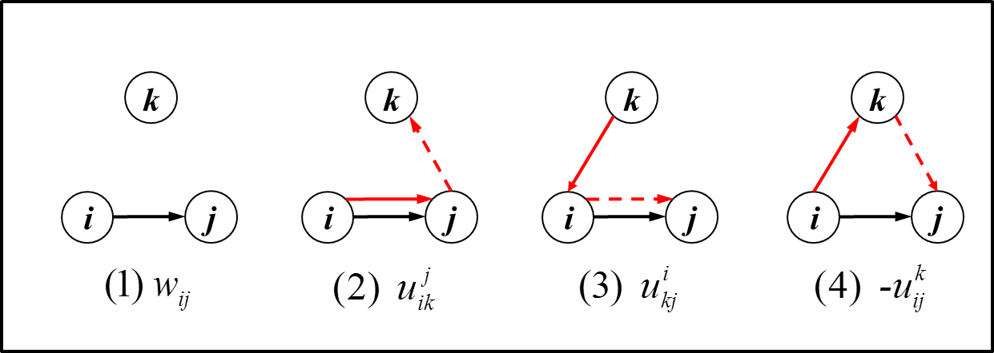}
   \end{center}
  \caption{Four Scenarios for Flow on Arc $(i, j)$ with Triples Variables}
  \label{fig:3}
\end{figure}

Figure \ref{fig:3} illustrates the calculation of the total flow on arc $(i,j)$, which can be derived from the triples variables as follows:

\begin{enumerate}
\item[(1)] If the delivery request from node $i$ to node $j$ is accepted, then that delivery puts $w_{ij}$ units of flow on arc $(i,j)$.
\item[(2)] If $u_{ik}^j$ is positive then the solution diverts $u_{ik}^j$ units of flow from node $i$ to node $k$ onto arc $(i,j)$.
\item[(3)] If $u_{kj}^i$ is positive, then the solution diverts $u_{kj}^i$ unit of flow from node $k$ to node $j$ onto arc $(i,j)$.
\item[(4)] If the solution diverts flow from node $i$ to node $j$ through node $k$, then $u_{ij}^k$ units of flow are diverted from arc $(i,j)$.
\end{enumerate}

\noindent Therefore, the total flow on arc $(i,j)$ is $w_{ij} y_{ij} + \sum_{(i,k,j) \in \mathcal{T}} u_{ik}^j + \sum_{(k,j,i) \in \mathcal{T}} u_{kj}^i - \sum_{(i,j,k) \in \mathcal{T}} u_{ij}^k$
where $\mathcal{T} = \{(i, j, k): i \in V \setminus \{n\}, j \in V \setminus \{1,i\}, k \in V \setminus \{1,n,i,j\} \}$.
\newpage
\noindent Example arc flows derived from the triples solution to the BPMP instance shown in Figure \ref{fig:2} are calculated in Table \ref{table:arcflows}.

\begin{table}[h!]
\centering
\begin{tabular}{|l|l|l|} \hline
Arc      & Expression for $\theta$  & Flow Value \\ \hline
$(1, 2)$ & $w_{12} y_{12} + u_{13}^2 + u_{14}^2 - u_{12}^3$  & $\theta_{12} = 0.4 + 0.3 + 0.2 - 0 = 0.9.$ \\  \hline
$(2, 3)$ & $w_{23} y_{23} + u_{24}^3 + u_{13}^2$             & $\theta_{23} = 0.3 + 0.3 + 0.3 = 0.9.$ \\ \hline
$(3, 4)$ & $w_{34} y_{34} + u_{14}^3 + u_{24}^3 - u_{34}^2 $ & $\theta_{34} = 0.4 + 0 + 0.3 - 0 = 0.7.$ \\ \hline
$(1, 3)$ & $w_{13} y_{13} + u_{12}^3 + u_{14}^3 - u_{13}^2$  & $\theta_{13} = 0.3 + 0 + 0 - 0.3 = 0.$ \\ \hline
$(3, 2)$ & $w_{32} y_{32} + u_{12}^3 + u_{34}^2$            & $\theta_{32} = 0 + 0 + 0 = 0.$ \\ \hline
\end{tabular}
\caption{Example Arc Flow Calculations}
\label{table:arcflows}
\end{table}

\subsection{BPMP MIP with Triples Variables}
\label{sec:TriplesMIP}
 The triples formulation of the BPMP replaces the $z$ variables of the node-arc formulation with triples variables. The multicommodity flow constraints  (\ref{eqn:2a})--(\ref{eqn:2i}) are replaced with the following set of {\em triples constraints} that relate the triples variables to arc flows:

\begin{equation}
\theta_{ij} = w_{ij} y_{ij} + \sum_{(i,k,j) \in \mathcal{T}} u_{ik}^j  + \sum_{(k,j,i) \in \mathcal{T}} u_{kj}^i -   \sum_{(i,j,k) \in \mathcal{T}} u_{ij}^k \qquad \forall (i, j) \in \mathcal{A}
\label{eqn:triples_constraints}
\end{equation}


\begin{eqnarray}
\theta_{ij} & \ge & 0  \qquad \qquad \forall (i, j) \in \mathcal{A} \label{eqn:theta_ge_0}  
\end{eqnarray}

 We impose the following constraints in order to force arc $(i,k)$ to be on the vehicle's route
 if variable $u_{ij}^k$ is positive:

\begin{equation}
u_{ij}^k  \le  Q x_{ik} \qquad \forall (i, j, k) \in \mathcal{T}  
\label{eqn:link_u_x}
\end{equation}

\noindent These constraints provide a logical linkage between the $u$ variables and the $x$ variables, and replace constraint set (\ref{eqn:2d}) of the node-arc formulation.   Finally, the triples variables must be nonnegative:

\begin{eqnarray}
u_{ij}^k & \ge & 0        \qquad \forall (i, j, k) \in \mathcal{T} \label{eqn:triples_ge_0}.  
\end{eqnarray}

\noindent The triples formulation maximizes profit (\ref{eqn:1}) subject to the vehicle-routing constraints, (\ref{eqn:2e})--(\ref{eqn:2h}), subtour elimination constraints, (\ref{eqn:2k}), variable domain constraints (\ref{eqn:s_ge_0})--(\ref{eqn:binary_y}), conditional arc flow constraints, (\ref{eqn:conditional arc flow}), and constraints (\ref{eqn:triples_constraints})--(\ref{eqn:triples_ge_0}) introduced above.
The complete expansion of the triples formulation for a generic four-node BPMP instance is given in Appendix \ref{appendix-sec2}.
From the expansion (and also from Table \ref{table:arcflows}) we can see that triples variable $u_{13}^2$ appears in exactly three of the triples constraints: it has coefficient $+1$ in the constraints for arcs $(1, 2)$ and $(2 ,3)$, and coefficient $-1$ in the constraint for arc $(1, 3)$.  Generalizing from this example we make the following important observation about the triples formulation.

\begin{obs}
\label{obs:3.1}
Triples variable $u_{ij}^k$ appears in exactly three of the triples constraints (\ref{eqn:triples_constraints}): it has coefficient $+1$ in the constraints for arcs $(i,k)$ and $(k,j)$, and coefficient $-1$ in the constraint for arc $(i,j)$.
\end{obs}

The triples formulation does not directly specify the relationship between the $x$ and $y$ variables in the way that constraints (\ref{eqn:2d})--(\ref{eqn:2b}) do in the node-arc formulation. However, that relationship can be deduced from the triples solution using the following results: Theorem \ref{thm:3.1}, which states that there is a logical relationship between the $x$ and $y$ variables in a triples solution such that all accepted delivery requests indicated by positive $y$ variables are on the route indicated by the positive $x$ variables, and Theorem \ref{thm:trp_subscrip_sequence} which states that if $y_{ij}$ is positive then the route visits node $i$ prior to visiting node $j$.

\begin{theorem}
\label{thm:3.1}
If $y_{ij} = 1$ in a feasible solution to the triples formulation, then nodes $i$ and $j$ are both on the vehicle's route from node 1 to node $n$.
\end{theorem}

\proof{Proof}
Define set $\mathcal{N}_r$ to be all nodes on the route, including nodes 1 and $n$. Consider
a node $i$ that is not on the route. Since $i \notin \mathcal{N}_r$, $x_{ij} = 0$ for all $j > 1$ (i.e., $\{j| (i,j) \in \mathcal{A}\}$). Furthermore, $\theta_{ij} = 0$ for every $j > 1$ to satisfy (\ref{eqn:conditional arc flow}). Thus, the triples constraints for arcs emanating from $i$ are

\begin{eqnarray*}
&  w_{i2} y_{i2} + \sum_{(i,k,2) \in \mathcal{T}} u_{ik}^2  + \sum_{(k,2,i) \in \mathcal{T}} u_{k2}^i -   \sum_{(i,2,k) \in \mathcal{T}} u_{i2}^k  = 0, & (i, 2) \\
&  w_{i3} y_{i3} + \sum_{(i,k,3) \in \mathcal{T}} u_{ik}^3  + \sum_{(k,3,i) \in \mathcal{T}} u_{k3}^i -   \sum_{(i,3,k) \in \mathcal{T}} u_{i3}^k  = 0, & (i, 3) \\
& \vdots & \\
 &  w_{ij} y_{ij} + \sum_{(i,k,j) \in \mathcal{T}} u_{ik}^j  + \sum_{(k,j,i) \in \mathcal{T}} u_{kj}^i -   \sum_{(i,j,k) \in \mathcal{T}} u_{ij}^k  = 0, & (i, j) \\
& \vdots & \\
&  w_{i,n-1} y_{i,n-1} + \sum_{(i,k,n-1) \in \mathcal{T}} u_{ik}^{n-1}  + \sum_{(k,n-1,i) \in \mathcal{T}} u_{k,n-1}^i -   \sum_{(i,n-1,k) \in \mathcal{T}} u_{i,n-1}^k  = 0, & (i, n-1) \\
 &  w_{in} y_{in} + \sum_{(i,k,n) \in \mathcal{T}} u_{ik}^{n}  + \sum_{(k,n,i) \in \mathcal{T}} u_{kn}^i -   \sum_{(i,n,k) \in \mathcal{T}} u_{in}^k  = 0, & (i, n).
\end{eqnarray*}

\noindent From Observation \ref{obs:3.1}, summing the triples constraints for all arcs emanating from node $i$ yields an implied constraint

\begin{eqnarray*}
& & w_{i2} y_{i2} + \ldots + w_{ij} y_{ij} + \ldots + w_{i,n-1} y_{i,n-1} + w_{in} y_{in} +\\
& & \sum_{(k,2,i) \in \mathcal{T}} u_{k2}^i + \ldots + \sum_{(k,j,i) \in \mathcal{T}} u_{kj}^i + \ldots + \sum_{(k,n-1,i) \in \mathcal{T}} u_{k,n-1}^i  + \sum_{(k,n,i) \in \mathcal{T}} u_{kn}^i = 0.
\end{eqnarray*}

The $w$ parameters, $y$ variables, and $u$ variables are all nonnegative. As defined in Section \ref{sec:2.1}, $w_{ij} > 0$ for $(i,j) \in \mathcal{R}$, so the corresponding $y$ variables have to have $y_{ij} = 0$ for all $j$ to satisfy the equation. When $(i,j) \in \mathcal{A} \setminus \mathcal{R}$, the value of  the corresponding $y$ variables doesn't affect the result since there is no delivery request from $i$ to $j$.
By a similar argument, it follows that if $j \notin \mathcal{N}_R$ then $y_{ij} = 0$ for all $i$. Thus, $y_{ij} = 0$ for any delivery request starting or ending at a node that is not on the route. Conversely, the starting and ending nodes of every accepted delivery request are on the vehicle's route.
\Halmos
\endproof

\begin{theorem}
\label{thm:trp_subscrip_sequence}
If a triples variable $u_{ij}^k > 0$ in an optimal solution to the triples formulation,  then nodes $i$, $j$, and $k$ are all on the vehicle's route and the vehicle
visits node $i$ prior to visiting node $k$, and visits node $k$ prior to visiting node $j$ (i.e., $s_i < s_k < s_j)$.
\end{theorem}

\noindent We defer the proof of Theorem \ref{thm:trp_subscrip_sequence} to Appendix \ref{sec:Validity of the Enhanced Triples Formulation}.

\begin{theorem}
\label{thm:triples valid}
The triples formulation is a valid model of the BPMP.
\end{theorem}

\proof{Proof}
The triples formulation shares the objective function (\ref{eqn:1}) and vehicle routing constraints (\ref{eqn:2e})--(\ref{eqn:2h}) and (\ref{eqn:2k}) with the node-arc formulation. Constraint sets (\ref{eqn:conditional arc flow}) and (\ref{eqn:triples_constraints}) ensure that the vehicle is never carrying more than $Q$ tons of cargo.  Theorems \ref{thm:3.1} and \ref{thm:trp_subscrip_sequence} ensure that for a feasible triples solution there is a logical connection between the vehicle's route and the accepted delivery requests.
\Halmos
\endproof

\subsection{Enhanced Triples Formulation}
\label{sec:enhanced triples}
An important corollary to Theorem \ref{thm:triples valid} is that the nonnegativity constraints  on the arc-flow variables,  (\ref{eqn:theta_ge_0}), and the
constraints explicitly linking the routing and triples variables, (\ref{eqn:link_u_x}), are not necessary for the triples formulation to be valid.  That is, the proofs of Theorems \ref{thm:3.1} and \ref{thm:trp_subscrip_sequence} do not rely on
constraints (\ref{eqn:theta_ge_0}) and (\ref{eqn:link_u_x}).  In Appendix \ref{sec:Validity of the Enhanced Triples Formulation}, we describe a straight-forward procedure to convert an optimal solution with negative arc flows to an equivalent one in which all the arc flows are nonnegative. In our preliminary testing we found that relaxing (\ref{eqn:theta_ge_0}) led to faster solution times.

Dropping the linking constraints, (\ref{eqn:link_u_x}),  makes the triples formulation considerably more compact as there are on the order of $n^3$ triples in ${\cal T}$.  This reduces solution time, however it complicates the interpretation of the triples variables. For example, removing the linking constraints allows for an alternative solution to the BPMP instance shown in Figure \ref{fig:2} in which the $x$ and $y$ variables have the same values as before, as do the arc flows, but the positive triples variables are $u_{13}^2 = 0.5$, $u_{14}^3 = 0.2$, and $u_{24}^3 = 0.1$ The interpretation of  $u_{14}^3$ is complicated by the fact that arc $(1, 3)$ is not on the vehicle's route. With the linking constraints removed, $u_{ij}^k$ no longer represents the amount flow from $i$ to $j$ that is routed on arc $(i, k)$ and then a path from $k$ from $j$; instead, we generalize the definition of $u^k_{ij}$ to say it represents the amount flow from $i$ to $j$ that is composed of flow from $i$ to $k$ adjoined to flow from $k$ to $j$. In interpreting the alternative solution to the instance in Figure \ref{fig:2},  $u_{14}^3 = 0.2$ indicates that 0.2 tons are sent on an unspecified path from node 1 to node 3, and then on an unspecified path from node 3 to node 4. The path from node 3 to node 4 is easily resolved as arc $(3, 4)$ itself since $x_{34} = 1$.  As in the example, the path from node 1 to node 3 is resolved by noting that $u_{13}^2$ is positive, and arcs $(1, 2)$ and $(2, 3)$ are on the vehicle's path.

In our preliminary study  \citep{BO} we found  that the solution times with the triples formulation can be further reduced by adding the node-degree constraints (\ref{eqn:node-degree}) from the node-arc formulation and the following sets of valid inequalities:

\begin{eqnarray}
\sum_{(i, j) \in {\cal A} } w_{ij} y_{ij} \le Q, \quad i \in {\cal N} \setminus \{n\} \label{eqn:C1 1} \\
\sum_{(i, j) \in {\cal A}} w_{ij} y_{ij} \le Q, \quad j \in {\cal N} \setminus \{1\} .\label{eqn:C1 2}
\end{eqnarray}

\noindent The above are valid inequalities that are satisfied by any feasible solution because the total weight of the delivery requests accepted from or to a given node cannot exceed the vehicle capacity. However, this condition is not necessarily enforced by solutions to the LP relaxation with fractional $y_{ij}$ variables.

Following  \cite{DLMTZ}, the subtour elimination constraints (\ref{eqn:2k}) can be strengthened by lifting them to
\begin{equation}
      s_i - s_j + (n - 1) x_{ij} + (n-3) x_{ji} \le  n - 2  \qquad \forall  i \in {\mathcal N} \setminus \{1,n\}, j \in {\mathcal N \setminus \{1, i, n\}}. \label{eqn:lifted mtz}
\end{equation}

\noindent In our preliminary testing \citep{BO}, we found that lifting the MTZ constraints was beneficial for solving the triples formulation (but not beneficial for solving the node-arc formulation). Likewise, we found it beneficial to bound the node-sequence variables so that

\begin{equation}
1 \le s_i \le n \qquad \forall i \in {\mathcal N} \setminus \{1\}. \label{eqn:mtz upper bound}
\end{equation}

Hereinafter, we refer to the MIP resulting from making the following changes to the triples formulation as the {\em enhanced triples formulation}:
(i) add constraints (\ref{eqn:node-degree}), (\ref{eqn:C1 1}) and  (\ref{eqn:C1 2}), (ii) replace constraints (\ref{eqn:2k}) and  (\ref{eqn:s_ge_0}) with constraints (\ref{eqn:lifted mtz}) and (\ref{eqn:mtz upper bound}),  and (iii) drop constraints (\ref{eqn:theta_ge_0}) and (\ref{eqn:link_u_x}).  Tables \ref{table:columns} and \ref{table:rows} give upper bounds on the number of variables and structural constraints in the enhanced node-arc and enhanced triples formulations, respectively, for worst-case instances in which $\mathcal{R} = \mathcal{A}$.
As shown in Tables \ref{table:columns} and \ref{table:rows}, the enhanced triples formulation reduces the number of binary variables and the number of constraints in the MIP by factors of $n^2$ and $n$, respectively. Note that we give derivations of   $|\mathcal{A}| = n^2
- 3n + 3$ and $|\mathcal{T}| = n^3-7n^2+17n-14$ in Appendix \ref{appendix-notation}.

In Section \ref{sec:results}, we demonstrate how this reduction in MIP size leads to significant improvement in solution time compared to the enhanced node-arc model. 
It is worthwhile to note that the most significant factor determining the size of the enhanced node-arc formulation is the number of delivery requests, $|\mathcal{R}|$, which determines the number of $z$ variables and multicommodity flow constraints (\ref{eqn:2a})--(\ref{eqn:2j}).
However, the size of the enhanced triples formulation is determined primarily by the number of nodes and is essentially independent of $|\mathcal{R}|$. Therefore, the enhanced node-arc formulation may actually be smaller if the number of delivery requests is small relative to $n$.

\begin{table}[hbt]
\centering
\begin{tabular}{|c|c|c|c|c|}  \hline
Number of               & \multicolumn{2}{|c|}{Enhanced Node-Arc Formulation}     & \multicolumn{2}{c|}{Enhanced Triples Formulation} \\ \hline
Continuous              & $\theta_{ij}$    & $n^2 - 3n + 3$ & $\theta_{ij}$ & $n^2 -3n + 3$ \\  \cline{2-5}
Variables               & $s_i$      & $n$       & $s_i$      &   $n$  \\  \cline{2-5}
                        &  \multicolumn{2}{c|}{}  & $u_{ij}^k$ & $n^3-7n^2+17n-14$ \\ \hline
Total                   & \multicolumn{2}{|c|}{$n^2 - 2n + 3$} & \multicolumn{2}{|c|}{$n^3-6n^2+15n-11$} \\ \hline
Binary                  & $x_{ij}$ and $y_{ij}$ & $2(n^2 - 3n + 3)$ & $x_{ij}$ and $y_{ij}$ & $2(n^2 - 3n + 3)$ \\ \cline{2-5}
Variables               & $z_{kl,ij}$ & $n^4 - 6 n^3 + 15 n^2 - 18 n + 9$ & \multicolumn{2}{c|}{} \\ \hline
Total                   & \multicolumn{2}{|c|}{$n^4 - 6 n^3 + 17 n^2 - 24 n + 15$}  & \multicolumn{2}{|c|}{$2n^2 - 6n + 6$} \\ \hline
\end{tabular}
\caption{Comparison of Variable Counts}
\label{table:columns}
\end{table}

\begin{table}[hbt]
\centering
\begin{tabular}{|c|c|c|c|c|}  \hline
Number of  Constraints             & \multicolumn{2}{|c|}{Enhanced Node-Arc Formulation}     & \multicolumn{2}{c|}{Enhanced Triples Formulation} \\ \hline
Routing                            & (\ref{eqn:2e})--(\ref{eqn:2h}) & $n + 1$                & (\ref{eqn:2e})-(\ref{eqn:node-degree}) & $2 n - 1$\\  \hline
MTZ                                & (\ref{eqn:2k})                 & {$n^2 - 3n + 3$}       & (\ref{eqn:lifted mtz}) &  $n^2 - 5n + 6$ \\ \hline
Multicommodity Flow  & (\ref{eqn:2a})--(\ref{eqn:2i}) & $n^3 - 2n^2 + 3$ & (\ref{eqn:triples_constraints}) & $n^2 - 3n + 3$ \\  \hline
Conditional Arc Flow   & (\ref{eqn:conditional arc flow})                & $n^2 - 3n + 3$ & (\ref{eqn:conditional arc flow}) & $n^2 - 3n + 3$ \\  \hline
Valid Inequalities & \multicolumn{2}{c|}{---} &   (\ref{eqn:C1 1})--(\ref{eqn:C1 2})&  $2 n - 2$\\ \hline
Total              & \multicolumn{2}{|c|}{$n^3 - 5n + 10$}  & \multicolumn{2}{|c|}{$3n^2 - 7n + 9$}\\ \hline
\end{tabular}
\caption{Comparison of Constraint Counts}
\label{table:rows}
\end{table}

\newpage
\subsection{Upper Bound on Profit}
\label{sec:upper bound on profit}
We conclude this section by using the enhanced triples formulation to derive an upper bound on the maximum profit for a BPMP instance. In Section \ref{sec:results} we use this expression to demonstrate the strength of the formulation.
We derive the bound in Theorem \ref{thm:3.6} with the help of the following lemma.

\begin{lemma}
\label{lemma:alt_obj}
The BPMP objective function (\ref{eqn:1}) is equivalent to
\[
(p - c)  \sum_{(i,j) \in \mathcal{A}} d_{ij} w_{ij} y_{ij} - c \sum_{(i,j,k) \in \mathcal{T}} (d_{ik} + d_{kj}  - d_{ij}) u_{ij}^k - cv \sum_{(i,j) \in \mathcal{A}} d_{ij} x_{ij}.
\]
\end{lemma}

\proof{Proof}
Define variable  $\mu_{ij} = \sum_{(i,k,j) \in \mathcal{T}} u_{ik}^j + \sum_{(k,j,i) \in \mathcal{T}} u_{kj}^i  -  \sum_{(i,j,k) \in \mathcal{T}} u_{ij}^k$ and note that

\begin{equation}
\theta_{ij} = w_{ij} y_{ij} + \mu_{ij}. \label{eqn:thm4}
\end{equation}

\noindent Substituting (\ref{eqn:thm4}) into (\ref{eqn:1}), the objective function can be rewritten as

\begin{equation}
(p - c)  \sum_{(i,j) \in \mathcal{A}} d_{ij} w_{ij} y_{ij} - c \sum_{(i,j) \in \mathcal{A}} d_{ij} \mu_{ij} - cv \sum_{(i,j) \in \mathcal{A}} d_{ij} x_{ij}. \label{eqn:11}
\end{equation}

\noindent Since $u_{ij}^k$ appears in the triples constraints for $(i,j)$, $(i,k)$ and $(k,j)$,

\begin{equation}
\sum_{(i,j) \in \mathcal{A}} d_{ij} \mu_{ij} = \sum_{(i,j,k) \in \mathcal{T}} (d_{ik} + d_{kj} - d_{ij}) u_{ij}^k. \label{eqn:12}
\end{equation}

\noindent Substituting (\ref{eqn:12}) into (\ref{eqn:11}) the objective function becomes

\[
(p - c)  \sum_{(i,j) \in \mathcal{A}} d_{ij} w_{ij} y_{ij} - c \sum_{(i,j,k) \in \mathcal{T}} (d_{ik} + d_{kj}  - d_{ij}) u_{ij}^k - cv \sum_{(i,j) \in \mathcal{A}} d_{ij} x_{ij}.
\]
\Halmos
\endproof

We now apply Lemma \ref{lemma:alt_obj} to derive an upper bound of profit in any BPMP instance.

\begin{theorem}
The maximum profit for a BPMP instance is at most \[(pQ-cQ-cv) D.\]
\label{thm:3.6}
\end{theorem}

\proof{Proof}
As in the proof of Lemma \ref{lemma:alt_obj}, define $\mu_{ij} = \sum_{(i,k,j) \in \mathcal{T}} u_{ik}^j + \sum_{(k,j,i) \in \mathcal{T}} u_{kj}^i  -  \sum_{(i,j,k) \in \mathcal{T}} u_{ij}^k$. Since $\theta_{ij} = w_{ij} y_{ij} + \mu_{ij}$,  constraints (\ref{eqn:conditional arc flow}) and  (\ref{eqn:triples_constraints}) imply
\[
y_{ij} w_{ij} \le Q x_{ij} - \mu_{ij} \qquad \forall (i, j) \in \mathcal{A}.
\]

Thus, the objective function becomes
\begin{eqnarray*}
p \sum_{(i,j) \in \mathcal{A}}  d_{ij} w_{ij} y_{ij} - c \sum_{(i,j) \in \mathcal{A}} d_{ij} (w_{ij} y_{ij} + \mu_{ij}) - cv \sum_{(i,j) \in \mathcal{A}} d_{ij} x_{ij} & = & \\
(p - c) \sum_{(i,j) \in \mathcal{A}}  d_{ij} w_{ij} y_{ij} - c \sum_{(i,j) \in \mathcal{A}}  d_{ij} \mu_{ij} - cv \sum_{(i,j) \in \mathcal{A}} d_{ij} x_{ij} & \le & \\
(p - c) \sum_{(i,j) \in \mathcal{A}}  d_{ij} (Q x_{ij} -  \mu_{ij}) - c \sum_{(i,j) \in \mathcal{A}}  d_{ij} \mu_{ij} - cv \sum_{(i,j) \in \mathcal{A}} d_{ij} x_{ij} & = & \\
(p Q - c Q - c v) \sum_{(i,j) \in \mathcal{A}} d_{ij} x_{ij} - p \sum_{(i,j) \in \mathcal{A}} d_{ij} \mu_{ij}.
\end{eqnarray*}

From the time/distance limit constraint ($\sum_{(i,j) \in \mathcal{A}} d_{ij} x_{ij} \le D$) we have

\[
(p Q - c Q - c v) \sum_{(i,j) \in \mathcal{A}} d_{ij} x_{ij} - p \sum_{(i,j) \in \mathcal{A}} d_{ij} \mu_{ij} \le (p Q - c Q - c v) D - p \sum_{(i,j) \in \mathcal{A}} d_{ij} \mu_{ij}.
\]

\noindent Therefore the profit is at most $(p Q - c Q - c v) D - p \sum_{(i,j) \in \mathcal{A}} d_{ij} \mu_{ij}$. Applying equation (\ref{eqn:12}) from the proof of Lemma \ref{lemma:alt_obj}, profit $\le (p Q - c Q - c v) D - p \sum_{(i,j,k) \in \mathcal{T}} (d_{ik} + d_{kj} - d_{ij}) u_{ij}^k$. In the Euclidean case distances between nodes satisfy the Triangle Inequality $d_{ik} + d_{kj} - d_{ij} \ge 0$.  By definition $u_{ij}^k \ge 0$, so $p \sum_{(i,j,k) \in \mathcal{T}} (d_{ik} + d_{kj} - d_{ij}) u_{ij}^k \ge 0$, and the optimal profit is at most $(p Q - c Q - c v) D$.
\Halmos
\endproof

%% file: section-4-heuristic.tex
\section{Restricted Triples Heuristic}
\label{sec:heuristic}
\begin{figure}[h!]
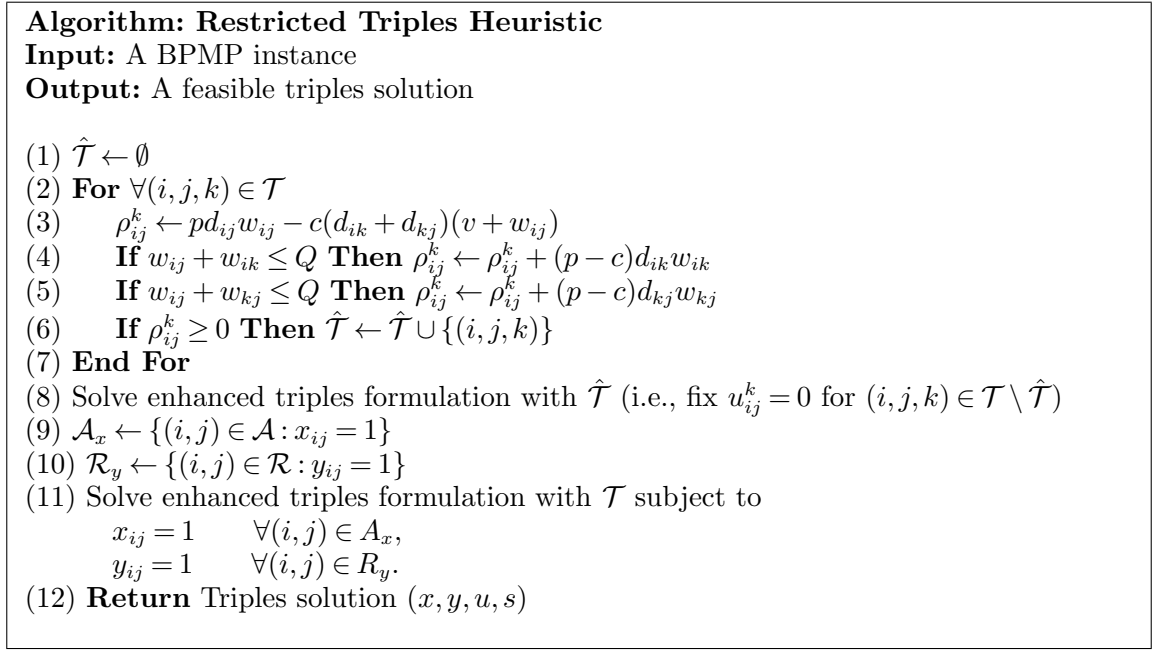

\anubox{
{\bf Algorithm: Restricted Triples Heuristic} \\
{\bf Input:}  A BPMP instance \\
{\bf Output:} A feasible triples solution\\
\\
(1) $\hat{{\mathcal T}} \leftarrow \emptyset$\\
(2) {\bf For} $\forall (i, j, k) \in \mathcal{T}$ \\
(3) \M $\rho^k_{ij} \leftarrow p d_{ij} w_{ij} - c(d_{ik} + d_{kj})(v + w_{ij})$ \\
(4) \M {\bf If} $w_{ij} + w_{ik} \le Q$ {\bf Then} $\rho^k_{ij} \leftarrow \rho^k_{ij}  + (p - c) d_{ik} w_{ik}$ \\
(5) \M {\bf If} $w_{ij} + w_{kj} \le Q$ {\bf Then} $\rho^k_{ij} \leftarrow \rho^k_{ij}  + (p - c) d_{kj} w_{kj}$ \\
(6) \M {\bf If} $\rho^k_{ij} \ge 0$ {\bf Then} $\hat{{\cal T}} \leftarrow \hat{{\cal T}} \cup \{(i,j,k)\}$ \\
(7) {\bf End For} \\
(8) Solve enhanced triples formulation with $\hat{{\cal T}}$ (i.e., fix $u^k_{ij} = 0$ for $(i,j,k) \in {\cal T} \setminus \hat{{\cal T}}$) \\
(9) ${\mathcal A}_x \leftarrow \{(i,j) \in {\mathcal A}: x_{ij} = 1\}$\\
(10) ${\mathcal R}_y \leftarrow \{(i,j) \in {\mathcal R}: y_{ij} = 1\}$\\
(11) Solve enhanced triples formulation with ${\cal T}$ subject to \\
\MM  $x_{ij} = 1 \qquad \forall (i,j) \in A_x$, \\
\MM  $y_{ij} = 1 \qquad \forall (i,j) \in R_y$.\\
(12) {\bf Return} Triples solution $(x, y, u, s)$ \\
}
\caption{Pseudo-code for Restricted Triples Heuristic}
\label{fig:heuristic}
\end{figure}

Our {\em restricted triples heuristic} for the BPMP solves the enhanced triples formulation with a relatively small subset of ``attractive" triples. It then fixes the positive $x$ and $y$ variables, and re-solves the MIP with the full set of triples to optimize the use of residual arc capacity on the fixed route. The heuristic is formalized in Figure \ref{fig:heuristic}.  The loop starting at line (2) and ending at line (7) calculates a {\em pseudo profit} $\rho^k_{ij}$ for each triple $(i,j,k)$. The initial value of $\rho^k_{ij}$ is the profit associated with accepting delivery request $(i, j)$ and transporting it on the two-arc path $i \rightarrow k \rightarrow j$. If the vehicle's capacity allows it, the pseudo profit is then incremented by the additional profit associated with accepting, and directly transporting, delivery requests $(i, k)$  and/or $(k, j)$.  The set of attractive triples, $\hat{\cal{T}}$, consists of the triples with nonnegative pseudo profit.

%% file: section-5-results.tex
\newpage
\section{Empirical Study And Analysis}
\label{sec:results}
In this section, we summarize our empirical study comparing solution times using the enhanced node-arc and triples formulations on BPMP instances from  \citep{Yu} with 10, 20, and 30 nodes, and new instances with 40 and 50 nodes from \citep{BPMPData}. In all instances, the delivery price, $p$, is \$1.20 per mile per ton, and the vehicle traveling cost, $c$, is \$1.00 per mile per ton. The average traveling speed of the vehicle is 50 miles per hour and the maximum time allowed for the backhaul trip, $\tau$, is 20 hours; thus, the maximum distance, $D$, the vehicle can travel is 1,000 miles. The capacity of the vehicle, $Q$, is 50 tons and the weight of the vehicle itself, $v$, is 5 tons. The delivery requests were generated randomly by taking 50 times a uniform random variable on the range [0,1] and rounding the result to one decimal place. The process for randomly determining the distances between nodes is described in  \citep{BO} and the data are available online \citep{BPMPData}.

\subsection{Comparison of MIP Sizes}
\label{sec:mip sizes}
Table \ref{table:mip sizes} details the relative MIP sizes, after reduction by AMPL and CPLEX's presolve routines, of the enhanced node-arc and  triples formulations for five example problems from our empirical study. As noted in Section \ref{sec:enhanced triples}, using the triples variables to represent the flow of delivery requests significantly reduces the number of binary variables and constraints in the MIP model. For example, the node-arc formulation of the 30-node instance in Table \ref{table:mip sizes} has over 600,000 binary variables and 25,000 constraints while the enhanced triples formulation has only 1,581 binary variables and 1,686 constraints. The enhanced triples formulation can have up to $n$ times more continuous variables than the enhanced node-arc formulation, however those variables don't contribute nearly as much to the MIP solution time as the binary variables.

\begin{table}[h]
\centering
\begin{tabular}{|c|c|r|c|r|} \hline
Nodes	&	\multicolumn{2}{c|}{Enhanced Node-Arc Formulation}	&	 		\multicolumn{2}{c|}{Enhanced Triples Formulation}			\\ \hline \hline
 	&	Continuous Variables	&	83	&	Continuous Variables	&	464	\\ \cline{2-5}
10	&	Binary Variables	&	4,377	&	Binary Variables	&	144	\\ \cline{2-5}
 	&	Constraints	&	832	&	Constraints	&	166	\\ \hline \hline
 	&	Continuous Variables	&	363	&	Continuous Variables	&	5,544	\\ \cline{2-5}
20	&	Binary Variables	&	106,306	&	Binary Variables	&	663	\\ \cline{2-5}
 	&	Constraints	&	7,262	&	Constraints	&	726	\\ \hline \hline
 	&	Continuous Variables	&	843	&	Continuous Variables	&	21,224	\\ \cline{2-5}
30	&	Binary Variables	&	617,836	&	Binary Variables	&	1,581	\\ \cline{2-5}
 	&	Constraints	&	25,292	&	Constraints	&	1,686	\\ \hline \hline
 	&	Continuous Variables	&	1,485	&	Continuous Variables	&	53,504	\\ \cline{2-5}
40	&	Binary Variables	&	2,040,418	&	Binary Variables	&	2,772	\\ \cline{2-5}
 	&	Constraints	&	60,846	&	Constraints	&	3,008	\\ \hline 
\end{tabular}
\caption{Example MIP Sizes}
\label{table:mip sizes}
\end{table}

\subsection{Comparison of Formulation Strength}
\label{sec:mip gaps}
Table \ref{table:mip gaps} compares the strength of the LP relaxations of the enhanced node-arc and triples formulations of 40 BPMP instances, ten instances for each value for $n \in \{10, 20, 30, 40\}$.  In Section \ref{sec:upper bound on profit}, we derived an upper bound on profit for a BPMP instance of at most $(pQ-cQ-cv) D$. Using the parameter values in our study, the upper bound is \$5,000. As shown in Table \ref{table:mip gaps}, the upper bound on profit from the LP relaxation of the enhanced triples formulation was consistently very close to \$5,000 while the bound from the enhanced node-arc formulation increased with the number of nodes ranging from an average of  \$12,642.50 for the smallest problem instances to  \$127,500.00 for the 40-node instances.  The gaps shown in Table \ref{table:mip gaps} are calculated relative to optimal MIP solutions. Thus, an instance with an optimal profit of \$3,550 and node-arc and triples upper bounds of \$13,000 and \$5,000, would have gaps of 266.20\% and 40.85\%, respectively. Table \ref{table:mip gaps} demonstrates that the enhanced triples formulation is stronger than enhanced node-arc formulation in addition to being more compact. 

\begin{table}[h]
\centering
\begin{tabular}{|c|c|r|r|c|r|} \cline{3-6}
\multicolumn{2}{c|}{} & \multicolumn{2}{c|}{Enhanced Node-Arc Formulation} & \multicolumn{2}{c|}{Enhanced Triples Formulation}\\ \hline
\multicolumn{2}{|c|}{Nodes}	 	&	LP Bound	&	\multicolumn{1}{c|}{Gap}	&	LP Bound	&	\multicolumn{1}{c|}{Gap}	\\ \hline \hline
	&	Min	&	\$5,000.00	&	13.5\%	&	\$5,000.00	&	13.5\%	\\ \cline{2-6}
10	&	Mean	&	\$12,642.50	&	192.43\%	&	\$5,000.00	&	47.12\%	\\ \cline{2-6}
 	&	Median	&	\$11,162.50	&	276.48\%	&	\$5,000.00	&	48.32\%	\\ \cline{2-6}
 	&	Max	&	\$24,600.00	&	969.89\%	&	\$5,000.00	&	152.61\%	\\ \hline \hline
	&	Min	&	\$35,000.00	&	620.58\%	&	\$5,000.00	&	2.94\%	\\ \cline{2-6}
20	&	Mean	&	\$47,500.00	&	1,003.10\%	&	\$5,001.05	&	14.43\%	\\ \cline{2-6}
 	&	Median	&	\$47,500.00	&	1,013.70\%	&	\$5,000.60	&	18.61\%	\\ \cline{2-6}
 	&	Max	&	\$60,000.00	&	1,337.33\%	&	\$5,003.00	&	23.98\%	\\ \hline \hline
	&	Min	&	\$57,500.00	&	1,287.57\%	&	\$5,000.00	&	5.57\%	\\ \cline{2-6}
30	&	Mean	&	\$80,000.00	&	1,671.72\%	&	\$5,002.75	&	14.45\%	\\ \cline{2-6}
 	&	Median	&	\$80,000.00	&	1,779.6\%	&	\$5,002.75	&	14.79\%	\\ \cline{2-6}
 	&	Max	&	\$112,500.00	&	2,273.99\%	&	\$5,005.03	&	21.93\%	\\ \hline \hline
	&	Min	&	\$92,500.00	&	2,087.34\%	&	\$5,000.05	&	18.07\%	\\ \cline{2-6}
40	&	Mean	&	\$127,500.00	&	3,114.21\%	&	\$5,000.08	&	27.29\%	\\ \cline{2-6}
 	&	Median	&	\$132,500.00	&	3,171.28\%	&	\$5,000.08	&	27.13\%	\\ \cline{2-6}
	&	Max	&	\$165,000.00	&	4,082.47\%	&	\$5,000.13	&	46.75\%	\\ \hline
\end{tabular}
\caption{Comparison of Strength of LP Relaxations}
\label{table:mip gaps}
\end{table}

\subsection{Comparison of Solution Times}
\label{sec:solution times}
We used AMPL version 10 to generate the MIPs, which were then solved with CPLEX version 12.6.0.0 on a Dell R730 computer with Dual 12 Core Intel Xeon@2.6GHz processors and 320GB of RAM. In \citep{BO} we found that a branching scheme that gives priority to the $x$ variables (routing decisions) over the $y$ variables (delivery-request decisions) improved solution time with the enhanced node-arc model compared to the default CPLEX settings. Other than this one change for the enhanced node-arc formulation, we used default settings for both AMPL and CPLEX.   We solved all of the 10-, 20-, 30-, and 40-node problems to optimality with both formulations. We solved all ten of the 50-node problems to optimality with the enhanced triples formulation, but due to excessive solution times only solved two of these instances to optimality with the enhanced node-arc model.

We report three measures of solution time: CPU time, real time, and “ticks”. By default, CPLEX version 12.6.0.0 uses a form of parallel processing that takes advantage of multiple cores and threads.   CPU time is the total time used by all threads on all processors of a CPLEX run, whereas real time (also known as wall clock time) is the amount of time that elapsed from the start of the run to the end of the run. Due to the nature of CPLEX's parallel processing and the fact that we ran our experiments on a multi-user system, we observed that multiple CPLEX runs with identical inputs showed variations in both CPU and real time. Therefore, we solved each problem instance three times with each formulation and compared the average CPU and real times reported. The ticks metric, also called “deterministic time”, is based on counting the number of instructions executed by the CPLEX solver and therefore shows no variation between multiple runs with the same inputs \citep{Ticks}. We report CPU time because it is a traditional performance measure in the literature, real time because it gives the most intuitive measure of computational effort, and ticks because it is a reproducible measure.

\subsubsection{CPU Time}
\label{sec:CPU}

Table \ref{table:CPU} summarizes and compares the average CPU times reported by CPLEX for the two formulations.
We can see that with its better upper bound and smaller constraint matrix, the enhanced triples formulation can indeed be solved much faster than the enhanced node-arc formulation. As reported in Table \ref{table:CPU}, the average solution times for the 10-node problem instances ranged from 0.78 to 5.95 seconds of CPU time using the enhanced node-arc formulation, and from 0.89 to 3.00 seconds of CPU time using the enhanced triples formulation. The {\em speedups}, the ratios of the average time using the enhanced node-arc formulation to the average time using the enhanced triples formulation, ranged from 0.39 to 4.52 with a median of 2.40 and geometric mean of 2.05. That is, on average, CPLEX solved the enhanced triples formulation approximately 2.05 times faster than the enhanced node-arc formulation for the 10-node instances. The average solution times for the 20-node problem instances ranged from approximately 126 seconds to 1,262 seconds (21 minutes) of CPU time using the enhanced node-arc formulation, and from 7.06 seconds to 50.76 seconds of CPU time using the enhanced triples formulation. The speedups for these instances, ranged from 6.33 to 137.11 with a median 68.98 and geometric mean of 43.90, respectively.  As expected, the average CPU time increased as a function of $n$ using both formulations. However, the rate of increase was much faster with the enhanced node-arc formulation. The median speedups for the 30-, and 40-node instances were  153.62, and 327.40, respectively, indicating that as the size of the problem instance grows it becomes increasingly faster to solve BPMP problem instances with the enhanced triples formulation.

\begin{table}[h!]
\centering
\begin{tabular}{|c|c|r|r|r|} \hline
\multicolumn{2}{|c|}{Nodes}	    	        &	Enhanced Node-Arc	&	Enhanced Triples	&	Speedup	\\ \hline \hline
	&	Min	                                &	0.78	            &	0.89	            &	0.39	\\ \cline{2-5}
10	&	Mean	                            &	4.36	            &	1.93	            &	2.05	\\ \cline{2-5}
 	&	Median	                            &	5.08	            &	1.91	            &	2.40	\\ \cline{2-5}
 	&	Max	                                &	5.95	            &	3.00	            &	4.50	\\ \hline \hline
	&	Min	                                &	126.00	            &	7.06	            &	6.33	\\ \cline{2-5}
20	&	Mean	                            &	725.00	            &	17.99	            &	43.90	\\ \cline{2-5}
 	&	Median	                            &	804.00	            &	11.43	            &	68.98	\\ \cline{2-5}
 	&	Max	                                &	1,262.00	        &	50.76	            &	137.06	\\ \hline \hline
	&	Min	                                &	20,454.00	        &	42.00	            &	34.33	\\ \cline{2-5}
30	&	Mean	                            &	55,459.00	        &	480.00	            &	139.02	\\ \cline{2-5}
 	&	Median	                            &	39,485.00	        &	413.00	            &	153.62	\\ \cline{2-5}
 	&	Max	                                &	140,354.00	        &	1,158.00	        &	872.08	\\ \hline \hline
	&	Min	                                &	334,025.00	        &	75.00	            &	44.46	\\ \cline{2-5}
40	&	Mean	                            &	6,643,337.00	    &	9,802.00	        &	480.58	\\ \cline{2-5}
 	&	Median	                            &	1,213,615.00	    &	5,162.00	        &	327.40	\\ \cline{2-5}
	&	Max	                                &	52,502,367.00	    &	50,649.00	        &	11,929.94	\\ \hline
\end{tabular}
\caption{Comparison of Average CPU Times (in seconds)}
\label{table:CPU}
\end{table}

\subsubsection{Real Time}
\label{sec:Real}

Table \ref{table:Real} summarizes and compares the average real times reported by CPLEX for the two formulations.  Measured in real time, which we argue is the most important metric to users, nearly all instances were solved faster with the enhanced triples formulation than the enhanced node-arc formulation.  The speedups were modest for the 10-node instances, but increased rapidly with problem size. For example,  the shortest average solution for the 40-node instances using the enhanced node-arc formulation was just over 5 hours while the longest average solution time using the enhanced triples formulation was approximately 30 minutes. Unless the delivery requests are known far in advance of the start of the vehicle's backhaul trip, the enhanced node-arc formulation is clearly impractical for instances with 40 or more nodes.

\begin{table}[h!]
\centering
\begin{tabular}{|c|c|r|r|r|} \hline
\multicolumn{2}{|c|}{Nodes}		 	&	Enhanced Node-Arc	&	Enhanced Triples	&	SpeedUp	\\ \hline \hline
	&	Min	                        &	0.21	            &	0.09	            &	0.86	\\ \cline{2-5}
10	&	Mean	                    &	0.65	            &	0.36	            &	2.12	\\ \cline{2-5}
 	&	Median	                    &	0.69	            &	0.28	            &	2.22	\\ \cline{2-5}
 	&	Max	                        &	0.96	            &	1.06	            &	4.00	\\ \hline \hline
	&	Min	                        &	31.00	            &	1.35	            &	7.98	\\ \cline{2-5}
20	&	Mean	                    &	61.00	            &	3.21	            &	18.90	\\ \cline{2-5}
 	&	Median	                    &	55.00	            &	3.02	            &	20.01	\\ \cline{2-5}
 	&	Max	                        &	105.00	            &	5.89	            &	46.20	\\ \hline \hline
	&	Min	                        &	1,594.00	        &	10.00	            &	30.96	\\ \cline{2-5}
30	&	Mean	                    &	3,198.00	        &	39.00	            &	87.09	\\ \cline{2-5}
 	&	Median	                    &	2,564.00	        &	37.00	            &	106.01	\\ \cline{2-5}
 	&	Max	                        &	6,166.00	        &	73.00	            &	212.53	\\ \hline \hline
	&	Min	                        &	18,413.00	        &	13.00	            &	56.61	\\ \cline{2-5}
40	&	Mean	                    &	329,773.00	        &	424.00	            &	353.30	\\ \cline{2-5}
 	&	Median	                    &	56,428.00	        &	265.00	            &	225.86	\\ \cline{2-5}
	&	Max	                        &	2,652,518.00	    &	1,863.00	        &	4,275.04	\\ \hline
\end{tabular}
\caption{Comparison of Average Real Times (in seconds)}
\label{table:Real}
\end{table}

\subsubsection{Ticks}
\label{sec:Ticks}

Table \ref{table:Ticks} summarizes and compares the number of ticks counted by CPLEX solving the problem instances in our data sets with the two MIP formulations. The ticks metric allows for a reproducible comparison of the computational effort required to solve BPMP instances for a given hardware configuration. Table  \ref{table:Ticks} shows that the effort increases with problem size at a much faster rate for the enhanced node-arc formulation than for the enhanced triples formulation. For example, the table shows a median speedup of 4.08, 25.39, 44.40, and 96.47, for the 10-, 20-, 30- and 40-node problem instances, respectively.

\begin{table}[h!]
\centering
\begin{tabular}{|c|c|r|r|r|} \hline
\multicolumn{2}{|c|}{Nodes}		 	    &	Enhanced Node-Arc	&	Enhanced Triples	&	Speedup	\\ \hline \hline
	&	Min	        &	41.91	        &	21.80	            &	1.47	\\ \cline{2-5}
10	&	Mean	    &	188.04	        &	45.47	            &	3.86	\\ \cline{2-5}
 	&	Median	    &	195.17	        &	41.09	            &	4.08	\\ \cline{2-5}
 	&	Max	        &	326.41	        &	84.09	            &	8.38	\\ \hline \hline
	&	Min	        &	18,475.00	        &	757.58	            &	8.89	\\ \cline{2-5}
20	&	Mean	    &	34,421.00	        &	1,573.38	        &	22.01	\\ \cline{2-5}
 	&	Median	    &	29,995.00	        &	1,316.83	        &	25.39	\\ \cline{2-5}
 	&	Max	        &	72,647.00	        &	3,614.30	        &	42.75	\\ \hline \hline
	&	Min	        &	812,789.00	        &	6,263.00	            &	16.06	\\ \cline{2-5}
30	&	Mean	    &	1,192,684.00	    &	31,266.00	            &	43.89	\\ \cline{2-5}
 	&	Median	    &	1,147,441.00	    &	28,218.00	            &	44.40	\\ \cline{2-5}
 	&	Max	        &	1,966,885.00	    &	58,062.00	            &	137.90	\\ \hline \hline
	&	Min	        &	6,601,682.00	    &	9,808.00	            &	34.79	\\ \cline{2-5}
40	&	Mean	    &	62,036,444.00	    &	257,478.00	            &	139.19	\\ \cline{2-5}
 	&	Median	    &	15,873,253.00	    &	177,890.00	            &	96.47	\\ \cline{2-5}
	&	Max	        &	463,811,772.00	    &	977,826.00	            &	1,417.66	\\ \hline
\end{tabular}
\caption{Comparison of Ticks}
\label{table:Ticks}
\end{table}

\subsection{Heuristic Results}

We ran the restricted triples heuristic described in Section \ref{sec:heuristic} on all of the problem instances in our data sets, and obtained optimal solutions for all but one case (one of the 50-node instances). The optimality gap of the one non-optimal solution was  1.58\% (\$4,427.55 vs. \$4,498.45). Table \ref{table:heuristic} gives summary statistics for the average real time to solve each instance three times with the heuristic, and speedups compared to exact solution using the full set of triples. As shown in the table, it turned out to be faster to solve some of the smaller problem instances with the exact approach, however the advantage of the heuristic became apparent when we solved the 40- and 50-node instances. The heuristic reduced the longest average solution time for the 40-node instances from approximately 30 minutes to a little under 9 minutes, and the maximum average time for the 50-node instances from approximately 83 minutes to just over 9 minutes.

\begin{table}[h!]
\centering
\begin{tabular}{|c|c|r|r|r|}  \cline{3-4}
\multicolumn{2}{c}{} & \multicolumn{2}{|c|}{Average Real Time (Seconds)} & \multicolumn{1}{c}{}\\ \hline
\multicolumn{2}{|c|}{Nodes}		 	    &	Enhanced Triples	&	Heuristic	        &	Speedup	\\  \hline \hline
	&	Min	                            &	0.09	            &	0.02	            &	0.77	\\ \cline{2-5}
10	&	Mean	                        &	0.36	            &	0.20	            &	1.68	\\ \cline{2-5}
 	&	Median	                        &	0.28	            &	0.19	            &	1.46	\\ \cline{2-5}
 	&	Max	                            &	1.06	            &	0.32	            &	6.63	\\ \hline \hline
	&	Min	                            &	1.35	            &	1.14	            &	0.62	\\ \cline{2-5}
20	&	Mean	                        &	3.21	            &	2.28	            &	1.43	\\ \cline{2-5}
 	&	Median	                        &	3.02	            &	2.02	            &	1.44	\\ \cline{2-5}
 	&	Max	                            &	5.89	            &	4.63	            &	2.86	\\ \hline \hline
	&	Min	                            &	10.00	            &	6.11	            &	1.24	\\ \cline{2-5}
30	&	Mean	                        &	39.00	                &	20.87	            &	1.97	\\ \cline{2-5}
 	&	Median	                        &	37.00	                &	17.18	            &	1.89	\\ \cline{2-5}
 	&	Max	                            &	73.00	                &	53.11	            &	3.74	\\ \hline \hline
	&	Min	                            &	13.00	                &	7.12	            &	1.79	\\ \cline{2-5}
40	&	Mean	                        &	424.00	                &	134.84	            &	2.78	\\ \cline{2-5}
 	&	Median	                        &	265.00	                &	104.96	            &	2.56	\\ \cline{2-5}
	&	Max	                            &	1,863.00	            &	517.44	            &	5.75	\\ \hline \hline
	&	Min	                            &	359.00	                &	152.02	            &	2.16	\\ \cline{2-5}
50	&	Mean	                        &	2,169.00	            &	256.69	            &	6.56	\\ \cline{2-5}
	&	Median	                        &	2,342.00	            &	214.00	            &	8.45	\\ \cline{2-5}
	&	Max	                            &	5,016.00	            &	558.23	            &	14.79	\\ \hline
\end{tabular}
\caption{Heuristic Speedups}
\label{table:heuristic}
\end{table}

%% file: section-6-conclusions.tex
\newpage
\section{Conclusions and Directions for Future Work}
\label{sec:conclusions}
We studied the problem of determining how to optimize profit on an empty delivery vehicle's backhaul trip to its depot, the backhaul profit maximization problem (BPMP).
We showed that the BPMP is ${\cal NP}$-hard.  The BPMP has previously been formulated in the literature as a mixed integer program based on the classical node-arc representation of multicommodity flow. We showed how the node-arc model can be strengthened and solved faster by employing preprocessing and valid inequalities, and in our computational experiments we found that solving the resulting enhanced node-arc formulation with the state-of-art MIP solver CPLEX is an effective, exact solution procedure for instances with up to 20 nodes and 343 delivery requests. To solve larger problem instances with up to 50 nodes and 2,353 delivery requests, we adapted a novel, compact representation of multicommodity flow to derive the triples and enhanced triples formulations of the BPMP. The enhanced triples formulation yields smaller MIP’s than the enhanced node-arc formulation, and appears from computational experiments to also have a stronger LP relaxation. In an empirical study, we demonstrated that  CPLEX can solve the enhanced triples formulation significantly faster than the enhanced node-arc formulation (e.g., approximately 350 times faster on average, in real time, for problem instances with 40 nodes and up to 1,483 delivery requests).

We also presented an easy-to-implement heuristic for the BPMP based on the enhanced triples formulation, and demonstrated that the heuristic can find optimal or near-optimal solutions to even the largest problems in our data set in less than 10 minutes of real time.   As a practical matter, we recommend the heuristic to 3PLs interested in solving large-scale problem instances. From our experience working with industry, we agree with \cite{Chandran2008} who argue that ``models and methodologies that can be easily implemented in a high-level modeling language are more likely to be implemented in practice, than specialized algorithms, that require sophisticated implementation."  Adopting the language from \cite{Chandran2008}, we recommend the heuristic as a way to ``solve the problem as efficiently as possible while working within the degrees of freedom offered by general purpose commercial solvers and modeling languages".

Finding exact solutions within practical solution-time limits to larger problem instances than those considered in our study, as well as solving generalizations of the BPMP such as considering multiple vehicles and/or individual time windows for the delivery requests will likely require developing specialized solution algorithms such as column-generation schemes or decomposition frameworks. Establishing the validity and computational superiority of the enhanced triples formulation in a straight-forward application of CPLEX (or comparable MIP solver) is a critical first step in these directions for future work.

%% file: appendix-A-notation.tex
\newpage
\section{Notation}
\label{appendix-notation}
\bigskip

\begin{tabular}{ll}
$n$ & number of nodes \\
$\mathcal{N} = \{1, 2, \ldots, n\}$ & set of nodes in the BPMP problem instance \\
$\mathcal{A}$ & set of arcs in the BPMP problem instance \\
$\mathcal{R}$ & set of delivery requests in the BPMP problem instance \\\\
$d_{ij}$ & distance in miles from node $i$ to node $j$ \\
$w_{ij}$ & weight in tons of delivery request from node $i$ to node $j$ \\
$p$ & revenue received in dollars per mile per ton \\
$v$ & weight of the vehicle in tons \\
$Q$ & carrying capacity of the vehicles in tons \\
$c$ & travel cost in dollars per mile per ton \\
$\tau$ & time in hours allowed for the vehicle to reach depot \\
$D$ & maximum distance in miles that the vehicle can travel \\ \\
$x_{ij}$ & binary variable equal to 1 if the vehicle travels on arc $(i, j)$ \\
$y_{kl}$ & binary variable equal to 1 if the delivery request $(k, l)$ is accepted \\
$z_{kl, ij}$ & binary variable equal to 1 if delivery request $(k, l)$ is routed on arc $(i, j)$ \\
$\theta_{ij}$ & variable indicating the total flow (tons carried by the vehicle) on arc $(i, j)$\\
$s_i$ & relative position (sequence number) of node $i$ in the vehicle's route \\ \\
$\mathcal{T}$ & set of node triples  $\{(i, j, k): i \in V \setminus \{n\}, j \in V \setminus \{1,i\}, k \in V \setminus \{1,n,i,j\} \}$  \\
$u_{ij}^{\ell}$ & flow (tons of cargo) on paths composed of a path from $(i, \ell)$ followed by a path from $\ell$ to $j$ \\
\end{tabular}

\bigskip
\noindent Considering the two distinct cases where $(i, j)$ is an arc in $\mathcal{A}$ as shown in the table below, we can see that $|\mathcal{A}| = n^2 -3 n + 3$. \\
\begin{center}
\begin{tabular}{cll}
Case & Description & Number of Triples \\ \hline
1 & $i = 1$, $j \in \mathcal{N} \setminus \{1\}$        & $(n - 1)$  \\
2 & $1 < i < n$, $j \in \mathcal{N} \setminus \{1, i\}$ &  $(n - 2)(n - 2)$  \\ \hline
\end{tabular}
\end{center}

\bigskip
\noindent Considering the four distinct cases where $(i, j, k)$ is a triple in $\mathcal{T}$ as shown in the table below, we can see that $|\mathcal{T}| = n^3-7n^2+17n-14$. \\
\begin{center}
\begin{tabular}{cll}
Case & Description & Number of Triples \\ \hline
1 & $i = 1$, $j = n$,  $k \in \mathcal{N} \setminus \{1, n\}$        & $(n - 2)$  \\
2 & $i = 1$, $1 < j < n$,  $k \in \mathcal{N} \setminus \{1, j, n\}$ &  $(n - 2)(n-3)$  \\
3 & $1 < i < n$, $j = n$, $k \in \mathcal{N} \setminus \{1, i, n\}$  & $(n - 2)(n -3)$  \\
4 & $1 < i < n$, $1 < j \neq  i < n$, $k \in \mathcal{N} \setminus    \{1, i, j, n\}$  &  $(n - 2)(n - 3)(n - 4)$ \\ \hline
\end{tabular}
\end{center}


%% file: appendix-B-mulitple-requests.tex
\section{BPMP Model with Multiple Delivery Requests Between Node Pairs}
\label{appendix-multiple-requests}

To allow for multiple requests between a given node pair, we can adopt conventional multicommodity flow notation whereby $\mathcal{R} =\{1, 2, 3, \ldots\}$ and each request $r \in \mathcal{R}$ has an origin, $o_r \in \mathcal{N}$, a destination (terminus) $t_r \in \mathcal{N}$, and a weight $w_r$.  Using this notation, the binary variable indicating whether or not a particular delivery request $r \in \mathcal{R}$ is accepted is $y_r$, and the objective function becomes

\[
p \sum_{r \in \mathcal{R}} d_{o_r, t_r} w_{r} y_{r} - c \sum_{(i, j) \in \mathcal{A}}  d_{ij} \theta_{ij} - c v \sum_{(i, j) \in \mathcal{A}} d_{ij} x_{ij}.
\]

\noindent In the node-arc formulation, the binary variable indicating whether or not request $r$ is transported on arc $(i, j)$ is $z_{r,ij}$ and constraints (\ref{eqn:2d})-(\ref{eqn:2i}) are rewritten as

\begin{eqnarray}
  \sum_{r \in \mathcal{R}} z_{r,ij} & \le & M x_{ij} \qquad \forall (i, j) \in \mathcal{A} \nonumber \\
  \sum_{j \in \mathcal{N} \setminus \{1, o_r\}} z_{r,o_r,j} &  = &  y_{r}   \qquad \forall r \in \mathcal{R} \nonumber \\
  \sum_{i \in \mathcal{N} \setminus \{t_r,n\}} z_{r,i,t_r} &  = &  y_{r} \qquad \forall r \in \mathcal{R} \nonumber\\
 \sum_{\{i \in \mathcal{N}: (i, h) \in \mathcal{A}\}} z_{r,ih} &  = & \sum_{\{j \in \mathcal{N}: (h, j) \in \mathcal{A}\}} z_{r,hj}  \qquad  \forall r \in \mathcal{R}, h \in \mathcal{N} \setminus \{o_r, t_r\} \nonumber \\
  \theta_{ij} & = & \sum_{r \in \mathcal{R}} w_{r} z_{r,ij} \qquad (i, j) \in \mathcal{A}. \nonumber
\end{eqnarray}

\noindent In the triples formulation, the triples constraints (\ref{eqn:triples_constraints}) are rewritten as

\begin{equation}
\theta_{ij} = \sum_{\{r \in \mathcal{R}: o_r = i, t_r = j\}} w_{r} y_{r} + \sum_{(i,k,j) \in \mathcal{T}} u_{ik}^j  + \sum_{(k,j,i) \in \mathcal{T}} u_{kj}^i -   \sum_{(i,j,k) \in \mathcal{T}} u_{ij}^k \qquad \forall (i, j) \in \mathcal{A}.
\nonumber
\end{equation}

\noindent Similar changes must also be made the enhanced node-arc and enhanced triples formulations. 

%% file: appendix-C-expanded-triples-model.tex
\section{The Enhanced Triples Formulation for a Generic Four-Node BPMP}
\label{appendix-sec2}

The objective function is

\begin{eqnarray*}
& & p (d_{12} w_{12} y_{12} + d_{13} w_{13} y_{13} + d_{14} w_{14} y_{14}  +  d_{23} w_{23} y_{23}  +  d_{24} w_{24} y_{24}  + d_{32} w_{32} y_{32}  +  d_{34} w_{34} y_{34}) -\\
& & c (d_{12} \theta_{12} + d_{13} \theta_{13} +  d_{14} \theta_{14}  +  d_{23} \theta_{23}  +  d_{24} \theta_{24}  + d_{32} \theta_{32}  +  d_{34} \theta_{34})  -  \\
& & cv (d_{12} x_{12} + d_{13} x_{13} + d_{14} x_{14}  +  d_{23} x_{23}  +  d_{24} x_{24}  + d_{32} x_{32}  +  d_{34} x_{34}).
\end{eqnarray*}

\bigskip

subject to {\bf routing constraints}:

\begin{eqnarray*}
 x_{12} + x_{13} + x_{14} & = & 1 \\
 x_{14} + x_{24} + x_{34} & = & 1 \\
 x_{12} - x_{23} - x_{24} + x_{32} & = & 0 \\
 x_{13} + x_{23} - x_{32} - x_{34} & = & 0 \\
 x_{12} + x_{32} &  \le & 1 \\
 x_{13} + x_{23} &  \le & 1
 \end{eqnarray*}

  {\bf subtour elimination constraints (lifted MTZ)}:

 \begin{eqnarray*}
 3x_{23} + s_2 - s_3 +x_{32}  &  \le &  2 \\
 3x_{32} + s_3 - s_2 +x_{23} &  \le &  2 \\
 \end{eqnarray*}

 {\bf distance/time limit constraint}
\begin{eqnarray*}
d_{12} x_{12} + d_{13} x_{13} + d_{14} x_{14} + d_{23} x_{23} + d_{24} x_{24} + d_{32} x_{32} + d_{34} x_{34} & \le &  D
\end{eqnarray*}

\smallskip

{\bf triples constraints}
\begin{eqnarray*}
\theta_{12} &  =  & w_{12} y_{12} + u_{13}^2 + u_{14}^2 - u_{12}^3 \\
\theta_{13} &  =  & w_{13} y_{13} + u_{12}^3 + u_{14}^3 - u_{13}^2 \\
\theta_{14} &  =  & w_{14} y_{14} - u_{14}^2 - u_{14}^3 \\
\theta_{23} &  =  & w_{23} y_{23} + u_{13}^2 + u_{24}^3 \\
\theta_{24} &  =  & w_{24} y_{24} + u_{14}^2 + u_{34}^2 - u_{24}^3 \\
\theta_{32} &  =  & w_{32} y_{32} + u_{12}^3 + u_{34}^2 \\
\theta_{34} &  =  & w_{34} y_{34} + u_{14}^3 + u_{24}^3 - u_{34}^2 \\
\end{eqnarray*}

\begin{eqnarray*}
u_{12}^3  & \ge & 0  \\
 u_{13}^2 & \ge & 0  \\
 u_{14}^2 & \ge & 0  \\
 u_{14}^3 & \ge & 0  \\
 u_{24}^3 & \ge & 0  \\
 u_{34}^2 & \ge & 0
\end{eqnarray*}

{\bf conditional arc flow constraints}

\begin{eqnarray*}
\theta_{12}  \le  Q x_{12} \\
\theta_{13}  \le  Q x_{13} \\
\theta_{14}  \le  Q x_{14} \\
\theta_{23}  \le  Q x_{23} \\
\theta_{24}  \le  Q x_{24} \\
\theta_{32}  \le  Q x_{32} \\
\theta_{34}  \le  Q x_{34}
\end{eqnarray*}

{\bf single node demand cuts}

\begin{eqnarray*}
 w_{12} y_{12} +  w_{13} y_{13} +  w_{14} y_{14} \le  Q  \\
  w_{23} y_{23} +  w_{24} y_{24} \le  Q  \\
  w_{32} y_{32} +  w_{34} y_{34} \le  Q  \\
  w_{12} y_{12} +  w_{32} y_{32} \le  Q  \\
  w_{13} y_{13} +  w_{23} y_{23} \le  Q  \\
      w_{14} y_{14} +  w_{24} y_{24} +  w_{34} y_{34} \le  Q
\end{eqnarray*}

{\bf MTZ upper bound}

\begin{eqnarray*}
   1 \le s_{2}  \le 3  \\
  1 \le s_{3}  \le 3  \\
  1 \le s_{4}  \le 3
\end{eqnarray*} 

%% file: appendix-D-proofs.tex
\newpage
\section{Validity of the Enhanced Triples Formulation}
\label{sec:Validity of the Enhanced Triples Formulation}

In this section, we complete the formal proof of the validity of the triples formulation and show that the enhanced triples formulation is also valid even though it relaxes the nonnegativity constraints  (\ref{eqn:theta_ge_0}) on the arc-flow variables. We begin in Section \ref{sec:DD} with a
description of a graphical representation of the positive $u$ variables in a triples solution called the {\em diversion digraph}. We derive several important properties of the diversion digraph that are then used in the proofs. In Section \ref{sec:proof theorem 3} we provide a proof of Theorem \ref{thm:trp_subscrip_sequence} from Section \ref{sec:TriplesMIP}, which states that vehicle route described by the $x$ variables visits nodes $i$, $j$, and $k$ in the correct order for every positive triples variable $u_{ij}^k$. We prove the validity of relaxing constraint set (\ref{eqn:theta_ge_0}) in the enhanced triples formulation in Section \ref{sec:free theta}.

\subsection{The Diversion Digraph}
\label{sec:DD}

\begin{figure}[h!]
    \begin{center}
   \includegraphics[scale=.4]{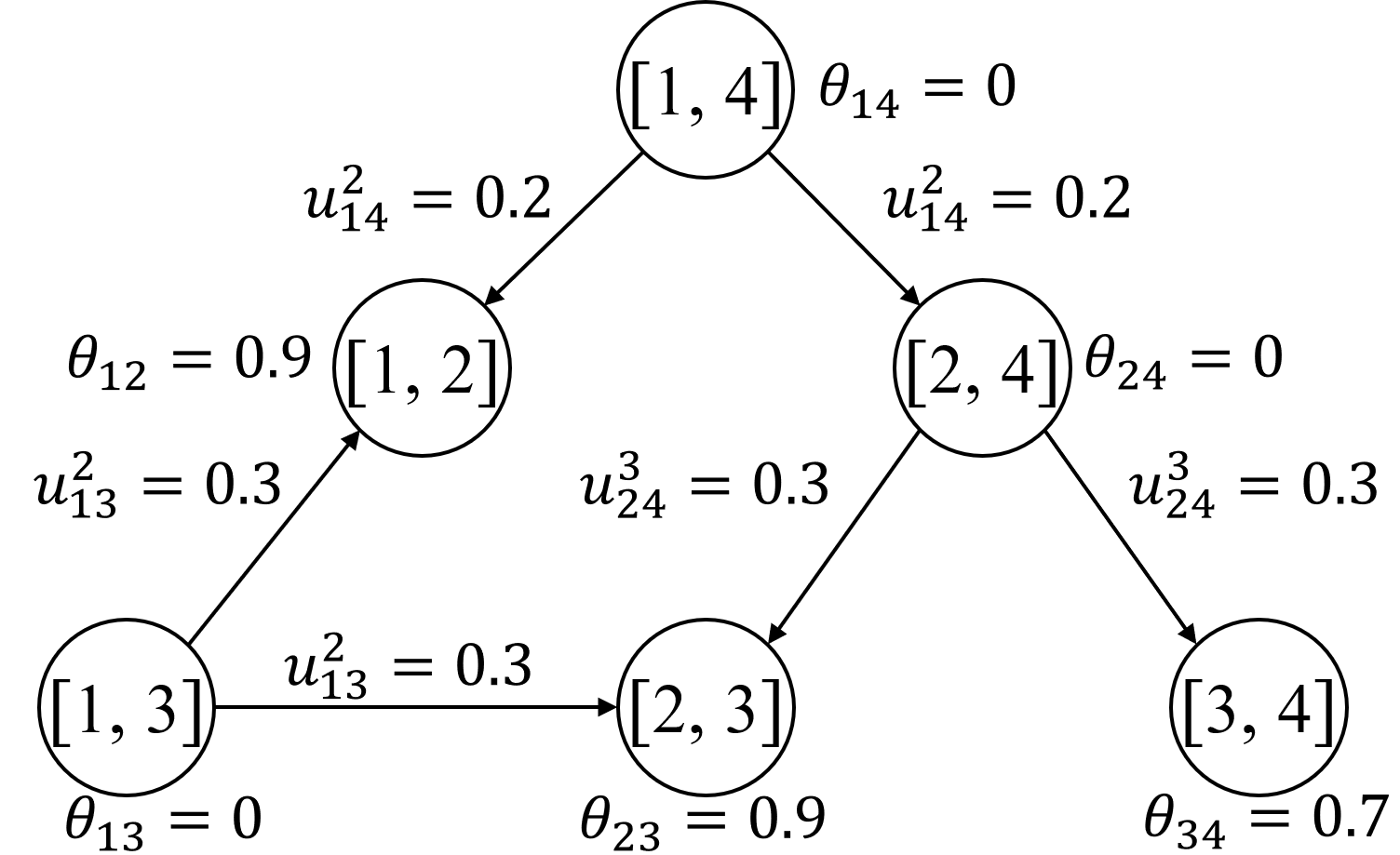}
   \end{center}
  \caption{Diversion Digraph for Triples Solution from Section \ref{sec:TriplesExample}}
  \label{fig:DD1}
\end{figure}

Consider a triples solution to a given BPMP problem. Arc $(i, j) \in {\cal A}$ is represented in the corresponding diversion digraph  by  a single node labeled $[i,j]$ if at least one of the triples variables on the right-hand side of the triples constraint  (\ref{eqn:triples_constraints}) for arc $(i, j)$ is positive.
The positive triples variable $u_{ij}^k$ is represented in the diversion digraph by two arcs emanating from the $[i, j]$ node: one to the $[i, k]$ node and another to the $[k, j]$ node. The diversion digraph for our example triples solution from Section \ref{sec:TriplesExample} is shown in Figure \ref{fig:DD1}. Since $u_{13}^2$ is positive, the diversion digraph has arcs from the $[1, 3]$ node to the $[1, 2]$ and $[2, 3]$ nodes. Likewise, the diversion digraph has arcs from the $[1, 4]$ node to the $[1, 2]$ and $[2, 4]$ nodes, and from the $[2, 4]$ node to the $[2, 3]$ and $[3, 4]$ nodes.

\begin{figure}[h!]
    \begin{center}
   \includegraphics[scale=.3]{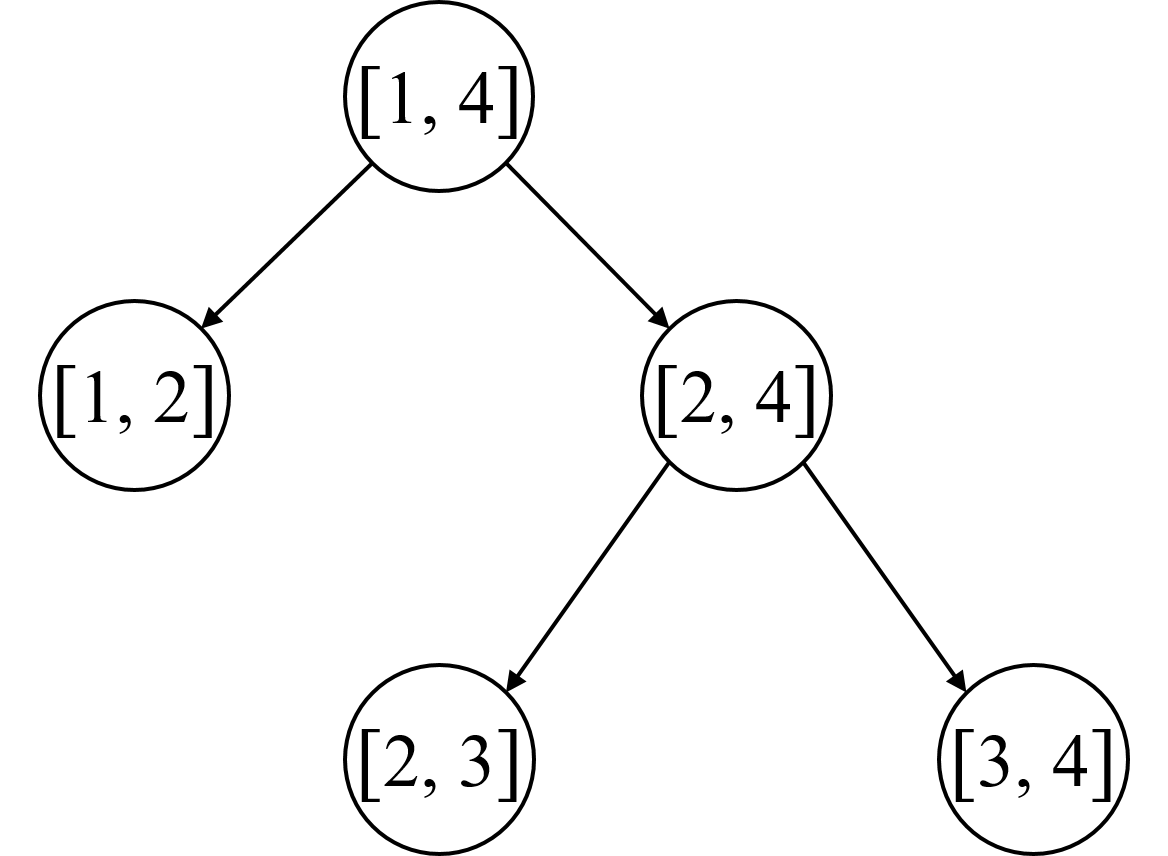}
   \end{center}
  \caption{Subtree of Diversion Digraph Rooted at (1, 4) Node}
  \label{fig:subtree}
\end{figure}

Nodes in the diversion digraph with out-degree zero are referred to as {\em leaf nodes}. For example, nodes $[1, 2]$, $[2, 3]$, and $[3, 4]$ in Figure \ref{fig:DD1} are leaf nodes. Observe that the leaf nodes in Figure \ref{fig:DD1} correspond precisely to arcs with positive flow in the triples solution. Furthermore, the leaf nodes in Figure \ref{fig:DD1} correspond to arcs on the vehicle's route from node $1$ to node $n$.  Recall from Section \ref{sec:TriplesExample} that the example triples solution accepts the delivery request from  node 1 to  node 4, and routes it on the path $1 \rightarrow 2 \rightarrow 3 \rightarrow 4$, which is represented by positive triples variables $u_{14}^2$ and $u_{24}^3$.  That is, flow from node 1 to node 4 is sent by the path $1 \rightarrow 2$ followed by a path from node 2 to node 4, and flow from node 2 to node 4 is sent on the two-arc path $2 \rightarrow 3 \rightarrow 4$.  This routing can be seen in the  diversion digraph in Figure \ref{fig:DD1} by considering the subtree of the diversion digraph rooted at the $[1, 4]$ node (shown in Figure \ref{fig:subtree}) and noting that the leaves of the subtree (i.e., $[1, 2]$, $[2, 3]$, and $[3, 4]$) correspond precisely to the arcs on the vehicle's route from node 1 to node 4. Likewise, the leaves of the subtree rooted at the $[1, 3]$ node correspond precisely to the arcs in the vehicle's path from node 1 to node 3: $(1, 2)$ and $(2, 3)$.

\begin{figure}[ht]
    \begin{center}
    \includegraphics[width=0.5\textwidth]{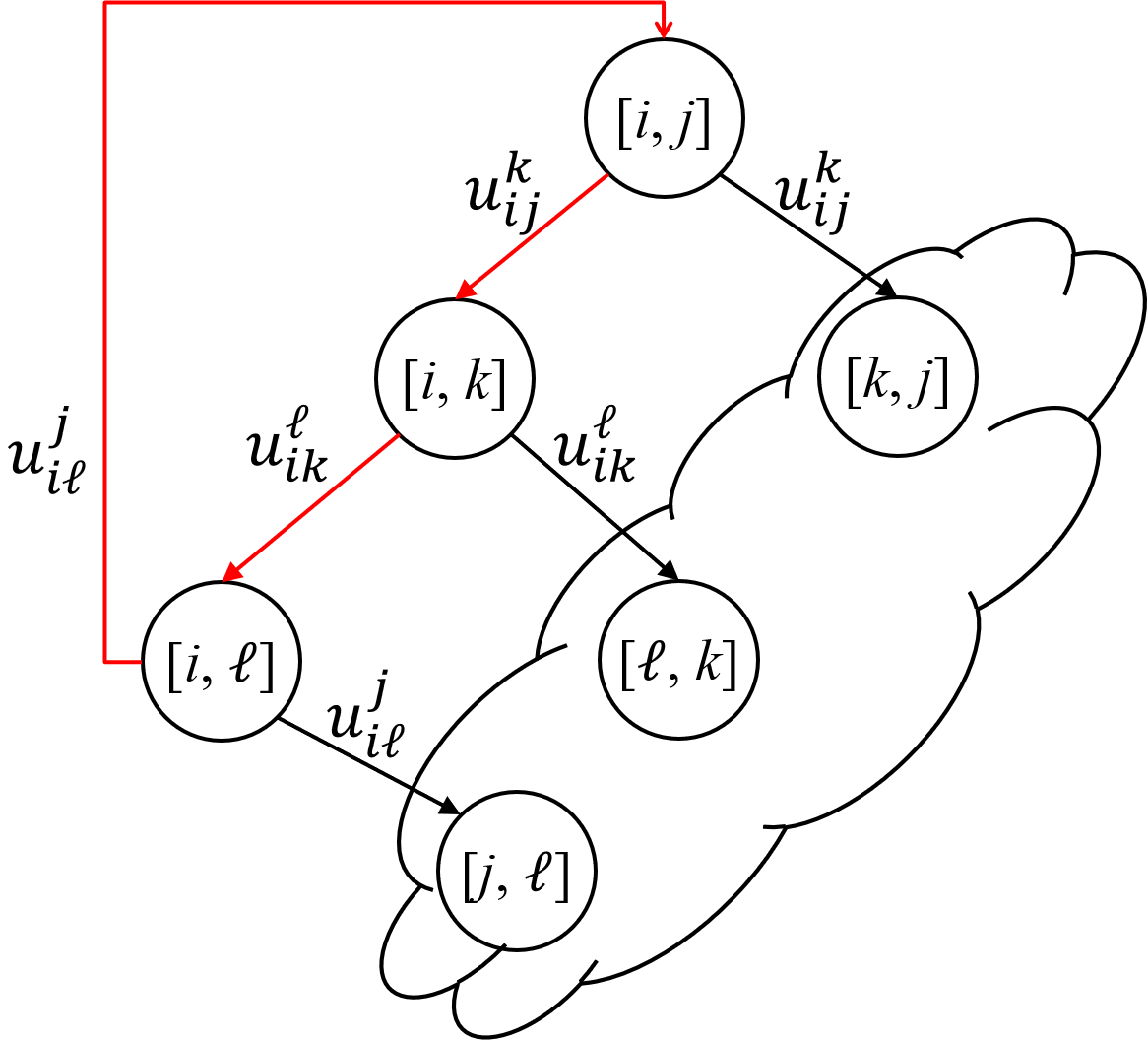}
   \end{center}
  \caption{Example Cycle in a Diversion Graph}
  \label{fig:cycle}
\end{figure}

An important feature of the diversion digraph is that it is acyclic if the corresponding triples solution is optimal. To see this,
suppose that the diversion digraph corresponding to an optimal solution to the enhanced triples formulation of a BPMP instance has at least one cycle as shown in Figure \ref{fig:cycle}. Note that the cloud shape around the $[k, j]$, $[\ell, k]$, and $[j, \ell]$, nodes indicates the rest of the subtree of the diversion digraph routed at the $[i, j]$ node. Let $\delta  = \min(u_{ij}^{k}, u_{ik}^{\ell}, u_{i\ell}^j) > 0$  denote the minimum  value of the triples variables whose corresponding nodes are on the cycle.  Decreasing each triples variable on the cycle by $\delta$ reduces one or more of the three triples variables to zero and breaks the cycle. Meanwhile, it follows from Observation \ref{obs:3.1} that the flow on arcs $(i, j)$, $(i, k)$, and $(i, \ell)$  in $\mathcal{A}$ are unchanged. That is, the affected triples variables appear as $(u_{i \ell}^j - u_{ij}^k)$, $(u_{ij}^k - u_{ik}^\ell)$,  and $(u_{ik}^{\ell}-u_{i \ell}^j)$ of the right-hand sides of the triples constraints (\ref{eqn:triples_constraints}) for arcs $(i,j)$, $(i, k)$, and $(i, \ell)$, respectively.  From  Observation \ref{obs:3.1}, the only other constraints involving variables $u_{ij}^{k}$, $u_{ik}^{\ell}$, $u_{i\ell}^j$ are the triples constraints for arcs $(k, j)$, $(\ell, k)$, and $(j, \ell)$, respectively. The values of $\theta_{kj}$, $\theta_{\ell k}$, and $\theta_{j \ell}$ can also be decreased by $\delta$ to preserve feasibility of the solution.  However, decreasing the value of the flow on arcs in $\mathcal{A}$ leads to an increase in the objective function value, which  contradicts the assumption that the initial triples solution is optimal. Generalizing from this example we state Theorem \ref{thm:acyclic_diversion_diagraph} without a formal proof.

\begin{theorem}
\label{thm:acyclic_diversion_diagraph}
The diversion digraph corresponding to an optimal solution to the enhanced triples formulation of the BPMP is acyclic.
\end{theorem}

\begin{theorem}
\label{thm:trp_subscrip_on_route}
If node $i$ is not on the vehicle's route in an optimal solution to the enhanced triples formulation of BPMP then $u^k_{ij} = 0$  for $\{j \in \mathcal{N}, k \in \mathcal{N}: (i,j,k) \in \mathcal{T}\}$, $u^i_{k j} = 0$  for $\{j \in \mathcal{N}, k \in \mathcal{N}: (k,j,i) \in \mathcal{T}\}$, and $u^{k}_{ji} = 0$  for $\{j \in N, k \in N: (j,i,k) \in \mathcal{T}\}$.
\end{theorem}

\proof{Proof}
Assume that there is a triple $(i, j, k) \in \mathcal{T}$ such that $u_{ij}^{k}  > 0$ and $i \notin \mathcal{N}_r$
where $\mathcal{N}_r$ is the set of nodes on the route as in the proof of Theorem \ref{thm:3.1}.

Since $i \notin \mathcal{N}_r$, $x_{i k}=0$, and $\theta_{i k} \le 0$ to satisfy the conditional arc-flow constraint (\ref{eqn:conditional arc flow}) for arc $(i, k)$. Thus, it follows that there must exist a triples variable $u_{i k}^{k_1} > 0$ to satisfy the triples constraint  (\ref{eqn:triples_constraints}) for $\theta_{i k}$.
Repeating the above analysis for $u_{i, k}^{k_1}$, we find another triples variable $u_{i, k_1}^{k_2} > 0$, and so forth. Let $\mathcal{T}_i = \{(i,j, k_0 = k), (i, k_0, k_1), \ldots, (i, k_h, k_{h + 1}), \ldots \} \subset \mathcal{T}$  be the sequence of triples identified by this process. Since
$\mathcal{T}$ is finite, the process must eventually reach an iteration $\ell + 1$ where the positive triple $(i, k_{\ell}, k_{\ell + 1})$ is already in $\mathcal{T}_i$ at which point $\cal{T}_i$ will contain a set of triples corresponding to a cycle in the diversion graph starting from node $[i,j]$ contradicting Theorem \ref{thm:acyclic_diversion_diagraph}.  Triples variables of the form $u_{k j}^i$ and $u_{j i}^k$ can be shown to be equal to zero by similar arguments.
\Halmos
\endproof

\subsection{Proof of Theorem \ref{thm:trp_subscrip_sequence}}
\label{sec:proof theorem 3}

Recall that Theorem \ref{thm:trp_subscrip_sequence} guarantees a logical connection between the routing ($x$) and triples ($u$) variables in the triples formulation by stating that

\begin{quote}
If a triples variable $u_{ij}^k > 0$ in an optimal solution to the triples formulation,  then nodes $i$, $j$, and $k$ are all on the vehicle's route and the vehicle
visits node $i$ prior to visiting node $k$, and visits node $k$ prior to visiting node $j$.
\end{quote}

\noindent Although the theorem was stated in the context of the triples formulation with nonnegative arc flows (constraint set (\ref{eqn:theta_ge_0})), it also holds for the enhanced triples formulation. 


\noindent {\it Proof of Theorem \ref{thm:trp_subscrip_sequence}}
Suppose that $u_{ij}^k > 0$ in an optimal triples solution. From Theorem \ref{thm:trp_subscrip_on_route} it follows that $i$, $j$, and $k$ are on the selected route. Now suppose further that the nodes are visited in the order $s_i < s_j < s_k$.  Since  $s_j < s_k$, it follows that $x_{kj}=0$ and $\theta_{kj} \le 0$. To satisfy the triples constraint (\ref{eqn:triples_constraints}) for arc $(k, j)$, there must be some node $\ell$
such that $u_{k j}^{\ell}>0$; and from Theorem \ref{thm:trp_subscrip_on_route}, $\ell$ is on the selected route. If $s_k < s_\ell$, then $s_j < s_\ell$ and $x_{\ell j} = 0$, and if $s_\ell < s_k$ then $x_{k \ell} = 0$. Hence, $\theta_{k \ell} \le 0$ or $\theta_{\ell j} \le 0$.  This process can be repeated until a sequence of positive triples variables is identified that corresponds to a cycle in the diversion graph contradicting Theorem \ref{thm:acyclic_diversion_diagraph}. A similar argument can be applied to the other four orderings of $i$, $j$, and $k$ that conflict with the definition of $u_{ij}^k$: $s_j < s_i < s_k$, $s_j < s_k < s_i$, $s_k < s_i < s_j$, and  $s_k < s_j < s_i$.
\Halmos

Theorem \ref{thm:trp_subscrip_sequence} establishes an important relationship between arcs with positive flow and leaf nodes in the diversion digraph.  This relationship is formalized in Theorem \ref{thm:pos_flow_leaf_node}  and  used to help justify relaxing the nonnegativity constraints on arc flow (\ref{eqn:theta_ge_0}) in Section \ref{sec:free theta}.

\begin{theorem}
\label{thm:pos_flow_leaf_node}
$\theta_{i j} > 0$ in an optimal solution to the enhanced triples formulation of BPMP if, and only if, node $[i, j]$ is a leaf node in the diversion digraph.
\end{theorem}

\proof{Proof} First, consider the if direction ($\theta_{ij} > 0 \Rightarrow [i,j]$ is a leaf node). Suppose that $\theta_{ij} > 0$ in an optimal triples solution, but node $[i, j]$ is not a leaf node in the corresponding diversion digraph. If  $[i, j]$ is not a leaf node then by definition there must be some $k \in \mathcal{N}_r$ such that $u_{ij}^k > 0$, and it follows from Theorem \ref{thm:trp_subscrip_sequence} that $s_i < s_k < s_j$. However, if $s_i <  s_k < s_j$ then $x_{ij} = 0$ violating the conditional arc-flow constraint (\ref{eqn:conditional arc flow}) for arc $(i, j)$.

Now consider the only if direction ($[i,j]$ is a leaf node $ \Rightarrow \theta_{ij} > 0$). In order for there to be a $[i, j]$ node in the diversion digraph,  there must be at least one positive triples variable of the form $u_{i k}^j$ or $u_{k j}^i$ in the corresponding triples solution. The fact that $[i, j]$ is a leaf node in the diversion digraph means that  $u_{i j}^k=0$ for all $(i, j, k) \in \mathcal{T}$.  Therefore, the right-hand side of the triples constraint  (\ref{eqn:triples_constraints}) for $(i, j)$ evaluates to a positive number.
\Halmos
\endproof

\subsection{Solutions with Negative Arc Flows}
\label{sec:free theta}

Our argument for the validity of relaxing the nonnegativity constraints on arc flow in the enhanced triples models has three main parts. In the first part we show that given a feasible solution in which $\theta_{ij} < 0$ for some arc $(i,j)$ there is a straight-forward way to derive another feasible solution with the same or better objective function value by increasing the flow on $(i, j)$ (Theorem \ref{lemma:2}).  Using this result, we show that  if $d_{ij} < d_{ik} + d_{kj}$ for every combination of $(i,j) \in \mathcal{A}$ and $\{k \in \mathcal{N}: \setminus \{i, j, k\} \in \mathcal{T}\}$ then all arc flows are nonnegative in any optimal solution to the enhanced triples formulation (Theorem \ref{thm:strict_triangle}). Theorem \ref{thm:strict_triangle} does not apply if there is an instance of a triangle equality in the driving distances such that  $d_{ij} = d_{ik} + d_{kj}$, however we show in this case that there exists an optimal solution in which all the arc flows are nonnegative. Furthermore, such a solution can be derived in a straight-forward manner from a solution with negative arc flows (Theorem \ref{thm:extended_triples_valid}).

\begin{theorem}
\label{lemma:2}
If $\theta_{ij}$ is negative for some arc $(i, j)$ in a feasible solution to the  triples formulation, then there must be at least one node $k$ such that the triples variable $u_{ij}^k$ is positive. Furthermore, the solution remains feasible if the values of $u_{ij}^k$, $\theta_{ik}$, and $\theta_{kj}$ are each reduced by $\delta = min(u_{ij}^k, -\theta_{ij})$, and the value of $\theta_{ij}$ is then increased by $\delta$.
\end{theorem}

\proof{Proof}
Consider the triples constraint (\ref{eqn:triples_constraints}) for arc $(i, j)$:

\[
\theta_{ij} = w_{ij} y_{ij} + \sum_{(i, k, j) \in \mathcal{T}} u_{ik}^j + \sum_{(k, j, i) \in \mathcal{T}} u_{k j}^i - \sum_{(i, j, k) \in \mathcal{T}} u_{i j}^k.
\]

\noindent Since the triples variables are nonnegative, it follows that if  $\theta_{ij} < 0$ then  $u_{ij}^k > 0$ for some $k$.  From Observation \ref{obs:3.1}, there are exactly three triples constraints that contain $u_{ij}^k$. If we decrease $u_{ij}^k$ by $\delta = min(u_{ij}^k, -\theta_{ij})$, then $\theta_{ij}$, $\theta_{ik}$, and $\theta_{kj}$ must be adjusted to satisfy the corresponding triples constraints.  Increasing $\theta_{ij}$ to  $\hat{\theta}_{ij} = \theta_{ij} + \delta$ preserves the equality of the left- and right-hand sides of the triples constraint for arc $(i, j)$.  Decreasing $\theta_{ik}$ and $\theta_{kj}$ to $\hat{\theta}_{ik} = \theta_{ik} - \delta$ and $\hat{\theta}_{kj} = \theta_{kj} - \delta$, respectively preserves equality for the triples constraints for arcs $(i, k)$ and $(k, j)$.

By construction, the reduced value of $u_{ij}^k$ is nonnegative. To complete the proof, we must show that the capacity constraints (\ref{eqn:conditional arc flow}) for arcs $(i,j)$, $(i, k)$, and $(k, j)$ remain satisfied after decreasing $u_{ij}^k$ by $\delta$.  This is clearly true for arcs $(i, k)$ and $(k, j)$ since $\hat{\theta}_{ik} < \theta_{ik} \le Q x_{ik}$
and $\hat{\theta}_{kj} < \theta_{kj} \le Q x_{kj}$; and since $\delta \le -\theta_{ij}$, the constraint holds for arc $(i,j)$ as $\hat{\theta}_{ij} \le 0 \le Q x_{ij}$.
\Halmos
\endproof

Note that since we assume that the driving distances in a BPMP instance satisfy the triangle inequality, it follows from Lemma \ref{lemma:alt_obj} that
applying Theorem \ref{lemma:2} increases the objective function value by $(d_{ik} + d_{kj} - d_{ij}) \delta \ge 0$.  It follows immediately that arc flows in an optimal solution to the enhanced triples formulation must be nonnegative for instances where $d_{ij}$ is always less than $d_{ik} + d_{kj}$ for all $k$.

\begin{theorem}
\label{thm:strict_triangle}
If $d_{ij} < d_{ik} + d_{kj}$ for every $(i, j, k) \in \mathcal{T}$ then all arc flows are nonnegative in any optimal solution to the enhanced triples formulation.
\end{theorem}

We conclude this section by outlining an algorithm to address the case of an optimal solution to the enhanced triples formulation in which $\theta_{ij} < 0$ for some $(i,j)$.  As noted above, this can only occur if $u_{ij}^k > 0$ for some node $k$ where $d_{ij} = d_{ik} + d_{kj}$. The algorithm  applies Theorem \ref{lemma:2} repeatedly until all the arc flows are nonnegative.

Suppose that $\theta_{ij} < 0$ for some $(i,j)$ in an optimal solution to the enhanced triples formulation. Rearranging the triples constraint for $(i, j)$ we have

\[
\sum_{(i, j, \ell) \in \mathcal{T}} u_{i j}^\ell - (-\theta_{ij}) = w_{ij} y_{ij} + \sum_{(i,\ell, j) \in \mathcal{T}} u_{i\ell}^j + \sum_{(\ell, j, i) \in \mathcal{T}} u_{\ell j}^i.
\]

\noindent Since $\theta_{ij} < 0$, the triples constraint implies that

\[
\sum_{(i, j, \ell) \in \mathcal{T}} u_{i j}^\ell \ge (-\theta_{ij}) > 0,
\]

\noindent which indicates that when $\theta_{ij} < 0$ there must exist a set of triples variables whose sum is greater than or equal to $-\theta_{ij}$.  Theorem \ref{lemma:2} can be applied to each of the  triples variables in the set to make the flow on arc $(i,j)$ increase to zero.
Since the given solution is optimal, it must be the case that $d_{i \ell} + d_{\ell j} - d_{ij} = 0$ for each  $u_{ij}^\ell > 0$ (Theorem \ref{thm:strict_triangle}); and so,  the profit is unchanged by each application of Theorem \ref{lemma:2}.

Thus, we have an alternative optimal solution in which $\theta_{ij} = 0$.  However, since $(i,j)$ is not a leaf node in the diversion digraph (Theorem \ref{thm:pos_flow_leaf_node}), we must check the new flow values on the arcs corresponding to the child nodes of node $[i, j]$ in the diversion digraph.

If child node $[i, k]$ of node $[i, j]$ is a leaf nodes then  $\theta_{ik} > 0$ prior to the application of Theorem \ref{lemma:2}.  Furthermore, since $\delta \le u_{ij}^k$, $\theta_{ik} \ge 0$ after the application and no adjustment is needed. However, if child node $[i, k]$ is not a leaf node, then $\theta_{ik} < 0$ after the application of Theorem \ref{lemma:2} and we must then apply Theorem \ref{lemma:2} to node $[i, k]$ and all of its immediate non-leaf node children, and so forth. Since the diversion digraph is acyclic (Theorem \ref{thm:acyclic_diversion_diagraph}), the adjustments starting from $[i, j]$ will eventually stop at leaf nodes. Thus, we can reduce the number of negative $\theta$ values in an optimal solution to the enhanced triples formulation by a finite number of applications of  Theorem \ref{lemma:2}.  Thus, we conclude with the following result establishing the validity of the enhanced triples formulation:

\begin{theorem}
\label{thm:extended_triples_valid}
There exists an optimal solution to the enhanced triples formulation in which $\theta_{ij} \ge 0$ for every arc $(i,j) \in \mathcal{A}.$
\end{theorem}